\documentclass[preprint,11pt,numbers]{elsarticle-arxiv}
\usepackage{lmodern}
\usepackage{amssymb,amsmath,amsthm,mathtools,mathrsfs}
\usepackage{enumitem}
\usepackage{indentfirst}
\usepackage{etoolbox}
\usepackage{url}
\usepackage[hidelinks]{hyperref}
\ifdefined\pdfgentounicode
  \InputIfFileExists{glyphtounicode.tex}{}{}
\fi
\setcitestyle{square,comma}
\numberwithin{equation}{section}
\newtheorem{theorem}{Theorem}[section]
\newtheorem{lemma}[theorem]{Lemma}
\newtheorem{proposition}[theorem]{Proposition}
\newtheorem{corollary}[theorem]{Corollary}

\theoremstyle{definition}
\newtheorem{definition}[theorem]{Definition}
\theoremstyle{remark}

\newcommand{\PP}{\mathbb P}
\newcommand{\EE}{\mathbb E}
\DeclareMathOperator{\Var}{Var}
\DeclareMathOperator{\Pois}{Pois}
\newcommand{\ind}{\mathbf 1}

\setlist[enumerate]{itemsep=3pt,topsep=5pt,parsep=0pt,listparindent=1.5em}
\setlist[itemize]{itemsep=3pt,topsep=5pt,parsep=0pt,listparindent=1.5em}
\patchcmd{\abstract}
  {\par\medskip\noindent\unskip\ignorespaces}
  {\par\medskip\ignorespaces}
  {}{\PackageError{manuscript}{Could not apply abstract indentation}{}}
\begin{document}
\bibliographystyle{alpha}
\begin{frontmatter}
\title{\textbf{\Large The Proof of the Alon--Kim Conjecture}}
\author[1]{Chengshuo Qiao\fnref{equalcontrib}}
\author[1]{Zhijun Lu\fnref{equalcontrib}}
\author[1]{Huimin Song\corref{cor1}}
\fntext[equalcontrib]{Chengshuo Qiao and Zhijun Lu contributed equally to this work.}
\cortext[cor1]{Corresponding author. \textit{E-mail address:} \url{hmsong@sdu.edu.cn} (H.Song).}
\affiliation[1]{organization={School of Mathematics and Statistics},
 addressline={Shandong University},city={Weihai},postcode={264209},
 state={Shandong},country={China}}
\begin{abstract}
In 1997, Alon and Kim proposed an asymptotic edge-coloring bound motivated
by a bound on the number of edges in intersecting uniform hypergraphs. The
Alon--Kim conjecture asserts that, for fixed integers $k\ge t\ge1$
and $\varepsilon>0$, every $k$-uniform multihypergraph $H$ in which any two
distinct edges intersect in at most $t$ vertices satisfies
$\chi'(H)\le(t-1+\frac{1}{t}+\varepsilon)\Delta(H)$ whenever $\Delta(H)$
is sufficiently large. In this paper, we prove the Alon--Kim conjecture. Our method
combines finite-dimensional convex optimization with probabilistic methods.
\end{abstract}
\end{frontmatter}
\section{Introduction}\label{s:intro}
\subsection{Preliminaries}\label{s:prelim}

All graphs, hypergraphs, and set families in this paper are finite. A
\emph{multihypergraph} $H=(V(H),E(H))$ has a vertex set $V(H)$ and an
edge multiset $E(H)$ consisting of nonempty subsets of $V(H)$.
Edges with the same vertex set are \emph{parallel}. Parallel edges are
distinguished by labels and regarded as distinct edges. All degrees and
edge counts are taken with multiplicity. We allow parallel edges unless
a hypothesis excludes them. We put
$d_H(v)=|\{e\in E(H):v\in e\}|$ and
$\Delta(H)=\max_{v\in V(H)}d_H(v)$, with $\Delta(H)=0$ when $V(H)$ is
empty. For distinct vertices $u,v$, their \emph{codegree} is
$d_H(u,v)=|\{e\in E(H):\{u,v\}\subseteq e\}|$. We write
$\Delta_2(H)=\max_{u\ne v}d_H(u,v)$, with value zero if $|V(H)|<2$.

The \emph{rank} of $H$ is $\max_{e\in E(H)}|e|$, with value zero when
$E(H)=\emptyset$. Thus $H$ has rank at most $k$ if $|e|\le k$ for every
$e\in E(H)$. The hypergraph $H$ is \emph{$k$-uniform} if $|e|=k$ for every
$e\in E(H)$, and \emph{$t$-simple} if $|e\cap f|\le t$ for every two
distinct edges $e,f\in E(H)$. The hypergraph $H$ is \emph{linear} if $H$ is $1$-simple,
and \emph{intersecting} if $e\cap f\ne\emptyset$ for every two distinct
edges $e,f\in E(H)$. A \emph{matching} is a subset
$I\subseteq E(H)$ whose distinct members have disjoint vertex sets. For
a positive integer $q$, a \emph{proper edge-coloring} is a map
$c:E(H)\to[q]$ such that $c(e)\ne c(f)$ whenever $e\ne f$ and
$e\cap f\ne\emptyset$, where $[q]=\{1,\ldots,q\}$. The least such $q$
is the \emph{chromatic index} $\chi'(H)$; we put $\chi'(H)=0$ if
$E(H)=\emptyset$.

A \emph{list assignment} gives a set $L(e)$ of allowed colors to each
edge $e\in E(H)$. A \emph{proper list edge-coloring} from the lists $L(e)$
is a map $c:E(H)\to\bigcup_{e\in E(H)}L(e)$ such that $c(e)\in L(e)$
for every $e\in E(H)$ and $c(e)\ne c(f)$ whenever $e\ne f$ and
$e\cap f\ne\emptyset$.

We also regard a family $\mathcal A$ of nonempty sets, with repetitions
allowed, as the edge multiset of a multihypergraph on $\bigcup\mathcal A$.
Repeated members are distinguished by labels, as for parallel edges.
The \emph{line graph} $U(\mathcal A)$ has vertex set $\mathcal A$,
and distinct $a,b\in\mathcal A$ are adjacent if $a\cap b\ne\emptyset$.
For a graph $U$, let $N_U[a]=\{a\}\cup\{b:ab\in E(U)\}$, and, for an
integer $R\ge0$, let
$B_U(a,R)=\{b:\operatorname{dist}_U(a,b)\le R\}$, where graph distance
is the minimum number of edges in a path and is infinite for vertices in
different components. For $\mathcal R\subseteq\mathcal A$, the graph
$U(\mathcal A)[\mathcal R]$ is the subgraph induced by $\mathcal R$.
A \emph{fractional matching} on $\mathcal A$ is a vector
$(x_a)_{a\in\mathcal A}$ satisfying $x_a\ge0$ for every $a$ and
$\sum_{a\ni v}x_a\le1$ for every underlying vertex $v$. Coordinates
with weight zero may be omitted. For any nonnegative vector $y$, the
\emph{support} is $\operatorname{supp}(y)=\{a:y_a>0\}$.

We write $\PP$, $\EE$, and $\Var$ for probability, expectation, and
variance. The indicator of an event $A$ is $\ind_A$. All conditional
probabilities $\PP(B\mid A)$ in this paper have $\PP(A)>0$ and mean
$\frac{\PP(A\cap B)}{\PP(A)}$. In every random experiment, the given hypergraph,
the weights on the hypergraph, and all numerical parameters selected as
functions of the hypergraph and the weights are fixed before the random
variables are sampled. Probabilities and expectations are taken only over
the explicitly specified random variables. A parameter may depend on the
entire given hypergraph without being a random variable.
For $\lambda\ge0$, the distribution $\Pois(\lambda)$ assigns mass
$\frac{e^{-\lambda}\lambda^j}{j!}$ to each integer $j\ge0$; the distribution
$\Pois(0)$ assigns probability one to zero. We put $(z)_+=\max\{z,0\}$ for real $z$.
All logarithms are natural unless a base is specified.

We write $X\sim\operatorname{Bern}(p)$ if $\PP(X=1)=p$ and
$\PP(X=0)=1-p$, where $0\le p\le1$. A uniform variable on $[0,1]$
has probability $v-u$ of belonging to $[u,v]$ for $0\le u\le v\le1$.
When several variables are attached to each member of a family, independence
between the groups means that the complete random vectors attached to
distinct members are mutually independent.

For events $(A_i)_{i\in I}$, a graph $\mathcal D$ with vertex set $I$
is a \emph{dependency graph} if, for each $i$, the event $A_i$ is
independent of the entire random vector
$(\ind_{A_j})_{j\notin N_{\mathcal D}[i]}$. If
each event is a function of a specified subset of mutually independent
random inputs, joining two indices whenever their events use a common
input gives such a graph.

For a bipartite graph $B$ with vertex classes $L,R$
and $S\subseteq L$, let $N_B(S)=\{y\in R:xy\in E(B)\text{ for some }x\in S\}$.

A \emph{system of distinct representatives} for an indexed family
$(S_c)_{c\in I}$ is an injective map $\phi:I\to\bigcup_{c\in I}S_c$
such that $\phi(c)\in S_c$ for every $c\in I$. For a $0$--$1$ matrix
$A$ with row set $R$ and column set $L$, the \emph{support} of column
$c$ is $S_c=\{i\in R:A_{i,c}=1\}$. If $L'\subseteq L$, a system
of distinct representatives for $(S_c)_{c\in L'}$ is equivalent to a
matching covering $L'$ in the bipartite graph with classes $L,R$ and
edges $ci$ whenever $A_{i,c}=1$.

For vectors $x,y\in\mathbb R^n$, write $x\cdot y=\sum_i x_i y_i$ and
$\|x\|_2=\sqrt{\sum_i x_i^2}$. A set $Q\subseteq\mathbb R^n$ is
\emph{convex} if $(1-\lambda)x+\lambda y\in Q$ for all $x,y\in Q$
and $0\le\lambda\le1$. A real function $F$ on a convex set $Q$ is \emph{concave}
if $F((1-\lambda)x+\lambda y)\ge(1-\lambda)F(x)+\lambda F(y)$
for all $x,y\in Q$ and $0\le\lambda\le1$; the function $F$ is
convex if $F((1-\lambda)x+\lambda y)\le(1-\lambda)F(x)+\lambda F(y)$
for all $x,y\in Q$ and $0\le\lambda\le1$.
The probability simplex is
$\Sigma_n=\{u\in\mathbb R^n:u_i\ge0\text{ for all }i,
\sum_i u_i=1\}$. If $Q\subseteq\mathbb R^n$ is a nonempty closed
convex set, the \emph{Euclidean projection} $\Pi_Q(z)$ is the unique
point of $Q$ minimizing $\|z-y\|_2$ over $y\in Q$. For
$Q=[0,2]^n$, the $i$th coordinate of $\Pi_Q(z)$ is $\min\{2,\max\{0,z_i\}\}$.
For a real matrix $A$, the Euclidean operator norm is
$\|A\|_2=\sup_{\|x\|_2=1}\|Ax\|_2$.

In the proof of \textbf{Lemma~\ref{lem:auxiliary}}, limits are taken as
$\Delta(H)\to\infty$ after $k,r,\rho,\xi$ and the constants chosen from
them have been fixed.

\subsection{Background and main results}\label{ss:background}

A proper edge-coloring partitions the edges of a hypergraph into matchings.
If a $k$-uniform multihypergraph $H$ has at least one edge, then
$\chi'(H)\ge\Delta(H)$. Each edge intersects at most
$k(\Delta(H)-1)$ other edges, so coloring the edges one at a time
gives $\chi'(H)\le k(\Delta(H)-1)+1$. For fixed $k$, the lower bound
$\Delta(H)$ and the upper bound $k(\Delta(H)-1)+1$ are linear in the
maximum degree. We seek the best leading coefficient and study how
restrictions on intersections between edges reduce the leading coefficient.

Frederickson, Hao, and Kelly~\cite{FredericksonHaoKelly2026} improved the
general upper bound as follows.

\begin{theorem}[\cite{FredericksonHaoKelly2026}]\label{thm:fhk}
For every integer $k\ge3$ and every $\varepsilon>0$, there is
$\Delta_0=\Delta_0(k,\varepsilon)>0$ such that every $k$-uniform
multihypergraph $H$ with $\Delta(H)\ge\Delta_0$ satisfies
\begin{equation}\label{eq:general-fhk-bound}
 \chi'(H)\le
 \left(k-\frac{k-1}{2k}
       +\frac{(k-1)(k-2)}{6k^2\sqrt{k}}
       +\varepsilon\right)\Delta(H).
\end{equation}
\end{theorem}

The coefficient in~\eqref{eq:general-fhk-bound}, excluding $\varepsilon$,
is smaller than $k$. The optimal coefficient is not known for every fixed
$k$.

To obtain a coefficient close to one, one can restrict the number of edges
containing a pair of vertices. A linear hypergraph satisfies
$\Delta_2(H)\le1$. More generally, it suffices to assume that
$\Delta_2(H)$ is a sufficiently small proportion of $\Delta(H)$.
In 1989, Pippenger and Spencer~\cite{PippengerSpencer1989} proved that
a sufficiently small codegree gives a coefficient arbitrarily close
to one for uniform hypergraphs whose vertex degrees are all close to the
maximum degree.

\begin{theorem}[\cite{PippengerSpencer1989}]\label{thm:pippenger-spencer}
For every integer $k\ge2$ and every $\alpha>0$, there exist $\delta>0$
and $\Delta_0>0$ such that every $k$-uniform multihypergraph $H$ satisfying
$\Delta(H)\ge\Delta_0$, $d_H(v)\ge(1-\delta)\Delta(H)$ for every
$v\in V(H)$, and $\Delta_2(H)\le\delta\Delta(H)$ has
$\chi'(H)\le(1+\alpha)\Delta(H)$.
\end{theorem}

Kahn~\cite{Kahn1996} extended \textbf{Theorem~\ref{thm:pippenger-spencer}}
to list edge-coloring and removed
the lower bound on vertex degrees. We use the following consequence for
proper edge-coloring; see also
\cite[Theorem~2.4]{KangKellyKuhnMethukuOsthus2023}.

\begin{theorem}[\cite{Kahn1996}]\label{thm:small-codegree}
For every positive integer $k$ and every $\alpha>0$, there exist $\xi>0$
and $\Delta_0>0$ such that every multihypergraph $J$ of rank at most $k$
with $\Delta(J)\ge\Delta_0$ and
$\Delta_2(J)\le\xi\Delta(J)$ satisfies
$\chi'(J)\le(1+\alpha)\Delta(J)$.
\end{theorem}

In particular, for fixed rank, linear multihypergraphs satisfy
$\chi'(H)=(1+o(1))\Delta(H)$ as $\Delta(H)\to\infty$.
More generally, \textbf{Theorem~\ref{thm:small-codegree}} gives
$\chi'(H)=(1+o(1))\Delta(H)$ whenever the rank is fixed and
$\Delta_2(H)=o(\Delta(H))$ as $\Delta(H)\to\infty$.
For $t\ge2$, however, $t$-simplicity does not imply the codegree hypothesis
of \textbf{Theorem~\ref{thm:small-codegree}}. For example, for a positive integer $m$, let $H$ consist
of the triples $\{u,v,w_j\}$, $j\in[m]$, where $u,v,w_1,\ldots,w_m$
are distinct. Then $H$ is $2$-simple and
$d_H(u,v)=\Delta(H)=m$.

In 1997, Alon and Kim~\cite{AlonKim1997} studied $k$-uniform,
$t$-simple hypergraphs without a small-codegree assumption. For an
intersecting hypergraph, every two distinct edges require
different colors, so $\chi'(H)=|E(H)|$. Alon and Kim proved the following bound
on the number of edges.

\begin{theorem}[\cite{AlonKim1997}]\label{thm:alon-kim-intersecting}
For every pair of integers $k\ge t\ge1$, there is $\Delta_0>0$ such that
every intersecting $k$-uniform, $t$-simple multihypergraph $H$ with
$\Delta(H)\ge\Delta_0$ satisfies
$|E(H)|\le(t-1+\frac{1}{t})\Delta(H)$.
\end{theorem}

Alon and Kim proposed the Alon--Kim conjecture in
\cite[Conjecture~1.1]{AlonKim1997}, asking whether the intersecting assumption
in \textbf{Theorem~\ref{thm:alon-kim-intersecting}} can be removed at the cost of
an arbitrarily small increase in the coefficient. More precisely, for fixed
integers $k\ge t\ge1$ and $\varepsilon>0$, the conjecture asserts that every
$k$-uniform, $t$-simple multihypergraph $H$ satisfies
$\chi'(H)\le(t-1+\frac{1}{t}+\varepsilon)\Delta(H)$ when $\Delta(H)$
is sufficiently large.

The case $t=1$ of the Alon--Kim conjecture follows from
\textbf{Theorem~\ref{thm:small-codegree}},
because every pair of vertices then belongs to at most one edge. For
$t\ge2$, the conjecture permits $\Delta_2(H)=\Delta(H)$. Recent progress includes
the bounds of Frederickson, Hao, and Kelly~\cite{FredericksonHaoKelly2026}:
if the maximum degree is sufficiently large, then
$\chi'(H)\le2.3581\Delta(H)$ for $3$-uniform, $2$-simple hypergraphs and
$\chi'(H)\le2.6791\Delta(H)$ for arbitrary $3$-uniform multihypergraphs.
The corresponding coefficients in the conjecture are $\frac{3}{2}$ and
$\frac{7}{3}$.

\par
In this paper, we prove the conjecture of Alon and Kim~\cite[Conjecture~1.1]{AlonKim1997}.

\begin{theorem}\label{thm:main}
For every pair of integers $k\ge t\ge1$ and every $\varepsilon>0$, there is $\Delta_0=\Delta_0(k,t,\varepsilon)>0$ such that every $k$-uniform, $t$-simple multihypergraph $H$ with $\Delta(H)\ge\Delta_0$ satisfies
\begin{equation}\label{eq:main-bound}
\chi'(H)\le \left(t-1+\frac{1}{t}+\varepsilon\right)\Delta(H).
\end{equation}
\end{theorem}

\subsection{Overview of the proof}\label{ss:outline}

For $t=1$, we apply \textbf{Theorem~\ref{thm:small-codegree}} directly.
Suppose that $r=t\ge2$, and put $C=r-1+\frac{1}{r}$ and
$\Delta=\Delta(H)$. We reduce the edge-coloring problem for $H$ to
edge-coloring problems for hypergraphs with small codegrees.
We partition $E(H)$ into families $E_1,\ldots,E_s$
and replace each family by a hypergraph $H_j$ of rank at most $k$.
We preserve intersections between edges in the same family so that we can
transfer a proper coloring of $H_j$ to the original edges in $E_j$.
We also bound the maximum degree of each $H_j$ in terms of
$\frac{\Delta}{s}$ so that using disjoint color sets for the $s$ families
still gives the required total number of colors. We prove the required
partition and degree bounds in \textbf{Lemma~\ref{lem:auxiliary}}.

\begin{lemma}\label{lem:auxiliary}
For every pair of integers $k\ge r\ge2$ and every $\rho,\xi>0$, there is
$\Delta_0=\Delta_0(k,r,\rho,\xi)>0$ such that the following holds.
Put $C=r-1+\frac{1}{r}$. Every $k$-uniform, $r$-simple multihypergraph
$H$ with $\Delta(H)\ge\Delta_0$ admits a partition
$E(H)=\dot\bigcup_{j=1}^s E_j$, multihypergraphs $H_1,\ldots,H_s$ of
rank at most $k$, and bijections $\varphi_j:E_j\to E(H_j)$ satisfying
\begin{enumerate}[label=\textnormal{(\roman*)}]
\item $\frac{\sqrt{\Delta(H)}}{2}\le s\le\sqrt{\Delta(H)}$;
\item $\Delta(H_j)\le\frac{(C+\rho)\Delta(H)}{s}$ and
$\Delta_2(H_j)\le\frac{\xi\Delta(H)}{s}$ for every $j\in[s]$;
\item if distinct $e,f\in E_j$ intersect, then
$\varphi_j(e)\cap\varphi_j(f)\ne\emptyset$.
\end{enumerate}
\end{lemma}

We choose $\xi$ from \textbf{Theorem~\ref{thm:small-codegree}} for a
small error $\alpha>0$. To apply \textbf{Theorem~\ref{thm:small-codegree}} at the degree scale in
\textbf{Lemma~\ref{lem:auxiliary}}, put
$q=\lceil\frac{(C+\rho)\Delta}{s}\rceil$. For each $j$, we add $q$
parallel singleton edges at a new vertex disjoint from $H_j$, obtaining a
multihypergraph $J_j$ with maximum degree exactly $q$ and the same
codegrees as $H_j$. By \textbf{Theorem~\ref{thm:small-codegree}}, $J_j$
has a proper edge-coloring with at most $(1+\alpha)q$ colors.
We restrict a proper edge-coloring of $J_j$ to $H_j$, transfer the colors
to $E_j$ using $\varphi_j$, and use disjoint color sets for different
indices $j$. The total number of colors is at most
$(1+\alpha)((C+\rho)\Delta+s)$. Since $s\le\sqrt{\Delta}$, we obtain
$(C+\varepsilon)\Delta$ colors by choosing $\alpha,\rho$ sufficiently
small and taking the maximum degree sufficiently large. We give the full
calculation in \textbf{Section~\ref{s:proof-main}}.

We now explain the construction in
\textbf{Lemma~\ref{lem:auxiliary}}. Choose
$0<\theta<\frac{1}{10(k+1)}$ according to~\eqref{eq:theta-choice},
and use $\eta=\theta^2$ in \textbf{Proposition~\ref{prop:rounding}}.
The construction has four parts.

\begin{enumerate}[label=\textnormal{(\roman*)}]
\item \textbf{Section~\ref{s:rounding}} constructs random matchings from
fractional matchings on multihypergraphs of rank at most $r$.
For positive weights $x_a$ with
$\sum_{a\ni v}x_a\le1$ at every vertex,
\textbf{Proposition~\ref{prop:rounding}} gives random matchings
$I_1,\ldots,I_P$ such that
\begin{equation}\label{eq:overview-marginals}
 \PP(a\in I_d)\le\frac{x_a}{C}\quad(d\in[P]),\qquad
 \frac{1}{P}\sum_{d=1}^P\PP(a\in I_d)
 \ge\frac{(1-\eta)x_a}{C}.
\end{equation}

Here $P$ depends only on $r,\eta$, and $d$ specifies a construction rather
than a sampled variable. We need both~\eqref{eq:overview-marginals} and a
restriction on the random inputs affecting each selection.
In \textbf{Definition~\ref{def:local-random-matching}}(iv), we require
equality of the selection indicators whenever two outcomes agree on all
labels and random variables in the specified retained neighborhood.
\textbf{Lemma~\ref{lem:rounding-resampling}} bounds the number and
location of decisions that change when one member's variables are resampled.
We apply \textbf{Lemma~\ref{lem:rounding-resampling}} to prove the
variance bound in \textbf{Lemma~\ref{lem:weighted-variance}}.

We first treat weights bounded below. We make several copies of each member
and impose vertex constraints within each copy and a constraint on the sum
over the copies of a member. Iterating the map in~\eqref{eq:gradient-iteration}
produces a vector whose total constraint violation is bounded in
\eqref{eq:constraint-excess}. Using~\eqref{eq:feasible-scaling} and the
coordinate deletions leading to~\eqref{eq:small-coordinate-gradient},
we obtain a feasible vector $y$. We use independent perturbations and
\textbf{Corollary~\ref{cor:hall-deficiency}} to bound, as
in~\eqref{eq:layer-deletion-bound}, the number of columns $c$ with $y_c>0$
that must be deleted so that the remaining columns can be
matched to distinct incident vertex rows.
We use the count in~\eqref{eq:layer-deletion-bound} to remove members in successive rounds
and then choose a matching by processing the removal rounds in reverse order. Poisson
sampling, as described in \textbf{Lemma~\ref{lem:poisson-expectation}},
removes the lower bound on positive weights. Finally, we apply
\textbf{Theorem~\ref{thm:brouwer}} to choose coefficients for which the lower bound
in~\eqref{eq:overview-marginals} holds for every member.

\item \textbf{Section~\ref{s:partition}} defines the sets that become
vertices of the auxiliary hypergraphs. We choose
$1=\delta_1>\cdots>\delta_{m_0+1}>0$, where
$m_0=1+\binom{k}{2}$. The level of an original edge $e$ is the least $i$
for which $\frac{d_H(u,v)}{\Delta}\notin[\delta_{i+1},\delta_i)$
for every pair of distinct vertices $u,v\in e$.
An index $i$ with the required property exists because the intervals are disjoint and $e$ contains fewer
than $m_0$ pairs of vertices. Let $\Gamma_i$ join vertices of codegree at
least $\delta_i\Delta$. We partition $e$ into the vertex sets $U$ of the
components of $\Gamma_i[e]$, and record each component vertex set $U$
as a part $a=(i,U)$. We write $U(a)=U$, $|a|=|U|$, and $B(e)$ for
the set of parts obtained from $e$, as in \textbf{Section~\ref{s:partition}}.

If distinct $a,b\in B(e)$ have level $i$, fewer than
$\delta_{i+1}\Delta$ original edges $f$ satisfy $a,b\in B(f)$.
If $d_a=|\{f\in E(H):a\in B(f)\}|$, then
$\sum_{a:v\in U(a)}\frac{d_a}{\Delta}\le1$ for every
original vertex $v$. We replace $\frac{d_a}{\Delta}$ by the positive weights $v_a$
in~\eqref{eq:part-weights}, preserving the sum bound and obtaining a positive
lower bound at each level. Thus the parts of size at most $r$ satisfy the
hypotheses of \textbf{Proposition~\ref{prop:rounding}}. A larger part has
$d_a=1$ by $r$-simplicity, and we select each part of size greater than $r$
using the Bernoulli variables defined in \textbf{Section~\ref{s:joint}}.

\item \textbf{Section~\ref{s:joint}} selects disjoint parts and retains
original edges whose parts are all selected. The marginal lower bounds
alone do not imply $\PP(B(f)\subseteq S_d)>0$: a random family could
choose exactly one of two disjoint parts, each with probability
$\frac{1}{2}$. We therefore sample independent Bernoulli variables $T_a$
with $\PP(T_a=1)=\theta v_a$. We include a part $a$ with $T_a=1$
only if $T_b=0$ for every other part $b$ intersecting $a$, and keep
$a\in I_d$ only if $T_b=0$ for every part $b$ intersecting $a$, including
$b=a$. The resulting
family $S_d$ is a matching, and~\eqref{eq:joint-probability-lower} gives
$q_{f,d}:=\PP(B(f)\subseteq S_d)\ge q_i>0$ for every edge $f$ of level $i$.

Let $\mathcal F_d$ consist of the edges $f$ such that
$\PP(a\in I_d)\ge\frac{(1-\theta)v_a}{C}$ for every $a\in B(f)$ with
$|a|\le r$. Each original edge belongs to $\mathcal F_d$ for at least
$(1-k\theta)P$ indices. For a part $a$, the sum $W_{a,d}$
in~\eqref{eq:weighted-sum} gives weight $\frac{1}{q_{f,d}}$ to each
$f\in\mathcal F_d$ using $a$ when all parts in $B(f)\setminus\{a\}$ are selected.
By testing selection only for the parts in $B(f)\setminus\{a\}$, we obtain the exact conditional mean
in~\eqref{eq:weighted-conditional-mean}.
\textbf{Lemma~\ref{lem:weighted-variance}} bounds the variance by
$\frac{K\delta_{i+1}\Delta^2}{q_i^2}$, with $K$ independent of all
codegree thresholds. We choose $\delta_{i+1}$ sufficiently small to obtain
\eqref{eq:conditional-overflow}. Consequently, conditional on
$B(f)\subseteq S_d$ for $f\in\mathcal F_d$, the probability that
$W_{a,d}>(A_\theta+2\theta)\Delta$ for some $a\in B(f)$ is at most
$\theta$. We retain $f\in\mathcal F_d$ only when $B(f)\subseteq S_d$ and
$W_{a,d}\le(A_\theta+2\theta)\Delta$ for every $a\in B(f)$.
The retention probability is bounded below in~\eqref{eq:edge-retention};
the upper bounds on $W_{a,d}$ will give the auxiliary degree bound
in~\eqref{eq:aux-degree-mean}.

\item \textbf{Section~\ref{s:auxiliary}} constructs the partition and
auxiliary hypergraphs. We take $s=P\lfloor\frac{\sqrt{\Delta}}{P}\rfloor$
independent experiments, with the same number for each index $d$.
Write $S^{(j)}$ for the selected matching in experiment $j$, whose index is
$d_j$. Let $R_{f,j}$ indicate retention of $f$ in experiment $j$, and put
$Z_f=\sum_j\frac{R_{f,j}}{q_{f,d_j}}$.
\textbf{Lemmas~\ref{lem:hoeffding} and~\ref{lem:lll}} give experiments
for which $Z_f\ge B_\theta s$ for every $f$, where
$B_\theta>0$ is defined in~\eqref{eq:theta-constants}.

We fix experiments satisfying $Z_f\ge B_\theta s$ for every $f$ and
assign each original edge independently
to experiment $j$ with probability $\frac{R_{f,j}}{q_{f,d_j}Z_f}$.
The auxiliary image of $f$ has vertex set $B(f)$. The retained weighted-sum bounds
give the expected degree bound in~\eqref{eq:aux-degree-mean}.
Item~\ref{part:pair} bounds the number of original edges using two parts
and gives the expected codegree bound in~\eqref{eq:aux-codegree-mean}.
We apply \textbf{Lemma~\ref{lem:bernstein}} and then
\textbf{Lemma~\ref{lem:lll}} to obtain assignments for which every
auxiliary degree is at most $\frac{(C+\rho)\Delta}{s}$ and every auxiliary
codegree is at most $\frac{\xi\Delta}{s}$. Finally, suppose that original edges
$f,g\in E_j$ meet at a vertex $u$. The parts of $B(f)$ and $B(g)$ containing
$u$ coincide because $B(f),B(g)\subseteq S^{(j)}$ and $S^{(j)}$ is a matching.
Hence $\varphi_j(f)\cap\varphi_j(g)\ne\emptyset$,
which proves \textbf{Lemma~\ref{lem:auxiliary}}(iii) and allows us to
transfer a proper coloring back to $H$.
\end{enumerate}

The remainder of this paper is organized as follows.
\textbf{Section~\ref{s:aux-lemmas}} states the auxiliary results.
\textbf{Section~\ref{s:rounding}} proves
\textbf{Proposition~\ref{prop:rounding}} and
\textbf{Lemma~\ref{lem:rounding-resampling}}.
\textbf{Sections~\ref{s:partition} and~\ref{s:joint}}
define and select the parts of original edges.
\textbf{Section~\ref{s:auxiliary}} proves
\textbf{Lemma~\ref{lem:auxiliary}}. Finally,
\textbf{Section~\ref{s:proof-main}} applies
\textbf{Theorem~\ref{thm:small-codegree}} to the auxiliary
hypergraphs and proves \textbf{Theorem~\ref{thm:main}}.

\section{Auxiliary Lemmas}\label{s:aux-lemmas}

We record the probability inequalities and the other auxiliary results used in the proof.

\begin{lemma}[\cite{Hoeffding1963}]\label{lem:hoeffding}
Let $X_1,\ldots,X_n$ be independent real random variables such that
$a_i\le X_i\le b_i$ almost surely, where $a_i,b_i$ are real numbers.
Put $S=\sum_{i=1}^nX_i$ and $V=\sum_{i=1}^n(b_i-a_i)^2$. If $V>0$,
then $\PP(S\le\EE S-u)\le \exp\left(-\frac{2u^2}{V}\right)$ for every $u>0$.
If $V=0$, then $S=\EE S$ almost surely, so
$\PP(S\le\EE S-u)=0$ for every $u>0$.
\end{lemma}

\begin{lemma}[\cite{BLM2013}]\label{lem:bernstein}
Let $X$ be a sum of finitely many independent Bernoulli random variables,
and put $\mu=\EE X$. For every $u>0$, we have
$\PP(X-\mu\ge u)\le\exp\left(-\frac{u^2}{2(\mu+\frac{u}{3})}\right)$.
For $X\sim\Pois(\mu)$ and every $u>0$, we also have
$\PP(X-\mu\ge u)\le\exp\left(-\frac{u^2}{2(\mu+\frac{u}{3})}\right)$.
\end{lemma}

We include a proof of the Poisson bound in \textbf{Lemma~\ref{lem:bernstein}}.
For $X\sim\Pois(\mu)$, we have
$\EE\exp(\lambda(X-\mu))=
\exp(\mu(e^\lambda-1-\lambda))$. For $0<\lambda<3$, the power
series gives $e^\lambda-1-\lambda\le
\frac{\lambda^2}{2(1-\frac{\lambda}{3})}$: indeed,
$j!\ge2\cdot3^{j-2}$ for every integer $j\ge2$. Applying Markov's
inequality to $\exp(\lambda(X-\mu))$, with
$\lambda=\frac{u}{\mu+\frac{u}{3}}$, gives
$\PP(X-\mu\ge u)\le\exp\left(-\frac{u^2}{2(\mu+\frac{u}{3})}\right)$ for
$\mu>0$. If $\mu=0$, then $X=0$ almost surely.

\begin{lemma}[\cite{BLM2013}]\label{lem:efron-stein}
Let $Y_1,\ldots,Y_n$ be independent random elements, and let
$X=f(Y_1,\ldots,Y_n)$ be a measurable real random variable with
$\EE X^2<\infty$. For each $i$, let $Y_i'$ be an independent copy of
$Y_i$, independent also of $Y_1,\ldots,Y_n$. Put
$Y^{(i)}=(Y_j^{(i)})_{j=1}^n$, where $Y_i^{(i)}=Y_i'$ and
$Y_j^{(i)}=Y_j$ for $j\ne i$, and let $X^{(i)}=f(Y^{(i)})$. Then
$\Var X\le\frac{1}{2}\sum_{i=1}^n\EE (X-X^{(i)})^2$.
\end{lemma}

\begin{lemma}[\cite{AlonSpencer2016}]\label{lem:lll}
Let $(A_i)_{i\in I}$ be finitely many events with a dependency graph of
maximum degree at most a nonnegative integer $b$. If $\PP(A_i)\le p$
for every $i\in I$ and $\mathrm e\,p(b+1)\le1$, then
$\PP(\bigcap_{i\in I}\overline{A_i})>0$.
\end{lemma}

We next state Hall's theorem and a bound on the number of vertices left
uncovered by a matching.

\begin{theorem}[\cite{Hall1935}]\label{thm:hall}
A finite bipartite graph $B$ with classes $L,R$ has a matching covering
every vertex of $L$ if and only if $|N_B(S)|\ge|S|$ for every
$S\subseteq L$.
\end{theorem}

\begin{corollary}\label{cor:hall-deficiency}
Let $B$ be a finite bipartite graph with classes $L,R$, and let $b\ge0$
be real. If $|N_B(S)|\ge |S|-b$ for every $S\subseteq L$, then there
are $L'\subseteq L$ and an injective map $\phi:L'\to R$ such that
$x\phi(x)\in E(B)$ for every $x\in L'$ and $|L\setminus L'|\le b$.
\end{corollary}

\begin{proof}
Put $h=\max_{S\subseteq L}(|S|-|N_B(S)|)$. The empty set shows
$h\ge0$, and the hypothesis gives $h\le b$. In particular, $h$ is a
nonnegative integer. Add $h$ new vertices to $R$, each adjacent to every
vertex of $L$, and denote the enlarged graph by $B^+$. For every nonempty
$S\subseteq L$, we have $|N_{B^+}(S)|=|N_B(S)|+h\ge|S|$; for
$S=\emptyset$, we have $|N_{B^+}(S)|=|S|=0$.
By \textbf{Theorem~\ref{thm:hall}}, $B^+$ has a matching $M^+$ covering $L$.
Remove from $M^+$ the edges incident to the new vertices, and let $M$
be the resulting matching in $B$. The edges of $M$ define $\phi$ and leave
at most $h\le b$ vertices of $L$ uncovered.
\end{proof}

We use~\eqref{eq:convex-minimizer-inequality} and
\eqref{eq:convex-projection-continuity} to estimate successive iterates in
\textbf{Section~\ref{s:rounding}} and to verify continuity when applying
\textbf{Theorem~\ref{thm:brouwer}}.

\begin{lemma}\label{lem:projection}
Let $Q\subseteq\mathbb R^n$ be nonempty, closed, and convex. For
$z,w\in\mathbb R^n$ and $p=\Pi_Q(z)$, we have
\begin{align}
 (z-p)\cdot(y-p)&\le0\qquad(y\in Q),\label{eq:convex-minimizer-inequality}\\
 \|\Pi_Q(z)-\Pi_Q(w)\|_2&\le\|z-w\|_2.\label{eq:convex-projection-continuity}
\end{align}
In particular, $\Pi_Q$ is continuous.
\end{lemma}

\begin{proof}
Choose $y_0\in Q$. The set
$\{y\in Q:\|z-y\|_2\le\|z-y_0\|_2\}$ is nonempty and compact.
The continuous function $y\mapsto\|z-y\|_2^2$ therefore attains a
minimum on $Q$. If distinct points $p,q\in Q$ both minimize the distance,
then $\frac{p+q}{2}\in Q$ and
\[
 \left\|z-\frac{p+q}{2}\right\|_2^2
 =\frac{\|z-p\|_2^2+\|z-q\|_2^2}{2}
   -\frac{\|p-q\|_2^2}{4}
 <\|z-p\|_2^2,
\]
contrary to minimality. Thus the minimizer $p=\Pi_Q(z)$ is unique.
For $y\in Q$ and $0\le t\le1$,
$p+t(y-p)\in Q$. The right derivative at $t=0$ of
$\|z-p-t(y-p)\|_2^2$ is therefore nonnegative, which gives
\eqref{eq:convex-minimizer-inequality}. Put $q=\Pi_Q(w)$. Applying
\eqref{eq:convex-minimizer-inequality} to $(z,p,q)$ and to $(w,q,p)$ gives
$(z-p)\cdot(q-p)\le0$ and $(w-q)\cdot(p-q)\le0$. Adding gives
$\|p-q\|_2^2\le(z-w)\cdot(p-q)\le\|z-w\|_2\|p-q\|_2$.
If $p=q$, both sides of \eqref{eq:convex-projection-continuity} are
nonnegative and the left side is zero. If $p\ne q$, dividing by
$\|p-q\|_2$ gives \eqref{eq:convex-projection-continuity}.
\end{proof}

\begin{theorem}[\cite{Brouwer1911}]\label{thm:brouwer}
Let $Q$ be a nonempty compact convex subset of a finite-dimensional real
Euclidean space. Every continuous map $T:Q\to Q$ has a point
$u\in Q$ satisfying $T(u)=u$.
\end{theorem}

\section{Constructing random matchings}\label{s:rounding}

Throughout \textbf{Section~\ref{s:rounding}}, $r\ge2$ is fixed and $C=r-1+\frac{1}{r}$.
We use the conventions for indexed families and line graphs from \textbf{Subsection~\ref{s:prelim}}. In particular, the distance between members $a,b$ is $\operatorname{dist}_{U(\mathcal A)}(a,b)$. The construction requires not only lower bounds on selection probabilities, but also a precise restriction on which random variables can affect each selection. We specify the random variables and selection functions in \textbf{Definition~\ref{def:local-random-matching}}.

The input family $\mathcal A$ and the fractional weights $x_a$ are fixed throughout each experiment. All constants and selection functions are chosen from $(\mathcal A,x)$ before any random variable is drawn and remain unchanged when a random variable is replaced by an independent copy. The experiment first uses Poisson variables to retain a subfamily of $\mathcal A$ and then chooses a matching within the retained subfamily.

\begin{definition}\label{def:local-random-matching}
Fix a finite indexed family $\mathcal A$ of nonempty sets, positive real numbers $x_a$ for $a\in\mathcal A$, and positive integers $m,M,R_*$. We consider random matchings $I$ obtained from the following random variables and measurable functions, subject to conditions~\textup{(i)}--\textup{(iv)}.
\begin{enumerate}[label=\textup{(\roman*)}]
\item Choose a positive integer $q_*$ and, for each $a\in\mathcal A$, a measurable function
$f_a:(\mathbb Z_{\ge0}\times[0,1]^{q_*})^{\mathcal A}\to\{0,1\}$.
The integer $q_*$ and the functions $f_a$ may depend on the fixed input $(\mathcal A,x)$ and on $m,M,R_*$, and are chosen before sampling. Draw $K_a\sim\operatorname{Pois}(mx_a)$ and $q_*$ uniform variables $U_{a,1},\ldots,U_{a,q_*}$ on $[0,1]$ for each $a$. The variables $K_a,U_{a,j}$, $a\in\mathcal A$, $j\in[q_*]$, are mutually independent. Write $K=(K_b)_{b\in\mathcal A}$, $U_a=(U_{a,1},\ldots,U_{a,q_*})$, and $\omega=((K_b,U_b):b\in\mathcal A)$ for the complete outcome.
\item Put $T_v=\sum_{b\ni v}K_b$ and retain precisely the family
\begin{equation}\label{eq:retained-family}
\mathcal A_K=\{a\in\mathcal A:K_a>0\text{ and }T_v\le M\text{ for every }v\in a\}.
\end{equation}
\item Put $X_a=f_a(\omega)$ and $I=\{a\in\mathcal A:X_a=1\}$. Require $X_a=0$ whenever $a\notin\mathcal A_K$, and require $X_aX_b=0$ whenever $a\ne b$ and $a\cap b\ne\emptyset$. These requirements say exactly that $I\subseteq\mathcal A_K$ and that $I$ is a matching.
\item For $a\in\mathcal A_K$, write
$B_K(a,R_*)=B_{U(\mathcal A_K)}(a,R_*)$.
For every $a\in\mathcal A$ and every two outcomes $\omega=((K_b,U_b))$ and $\omega'=((K_b',U_b'))$, require $f_a(\omega)=f_a(\omega')$ whenever $a\in\mathcal A_K\cap\mathcal A_{K'}$, $B_K(a,R_*)=B_{K'}(a,R_*)$, and $(K_b,U_b)=(K_b',U_b')$ for every $b\in B_K(a,R_*)$.
\end{enumerate}

In part~\textup{(iv)}, the labels and underlying sets belong to the fixed input. The equality $B_K(a,R_*)=B_{K'}(a,R_*)$ therefore also identifies the induced edges and input weights on the two balls. Part~\textup{(iv)} states precisely how the selection of $a$ is restricted to the retained ball $B_K(a,R_*)$. Retention in~\eqref{eq:retained-family} may depend on all incident Poisson variables through the sums $T_v$, including variables of members subsequently discarded. Once $\mathcal A_K$ is determined, the indicator $X_a$ must satisfy part~\textup{(iv)}. No constants or functions are chosen in response to the sampled outcome.
\end{definition}

\begin{proposition}\label{prop:rounding}
For every integer $r\ge2$ and every $0<\eta<\frac{1}{4}$, there are positive integers $P,m,M,R_*$ with the following property.
Let $\mathcal A$ be a finite indexed family of nonempty sets of size at most $r$, and suppose that $x_a>0$ and $\sum_{a\ni v}x_a\le1$ for every vertex $v$.
For this fixed input, there are $P$ collections of functions as in \textbf{Definition~\ref{def:local-random-matching}}, giving random matchings $I_d$ for $d\in[P]$, such that, writing $p_{a,d}=\mathbb P(a\in I_d)$,
\begin{equation}\label{eq:rounding-marginals}
0\le p_{a,d}\le\frac{x_a}{C}\quad(a\in\mathcal A,\ d\in[P]),
\qquad
\frac{1}{P}\sum_{d=1}^Pp_{a,d}\ge\frac{(1-\eta)x_a}{C}\quad(a\in\mathcal A).
\end{equation}

All four integers depend only on $r$ and $\eta$.
The functions may depend on the entire finite input $(\mathcal A,x)$, but are chosen before the random variables in \textbf{Definition~\ref{def:local-random-matching}} are drawn. The number $q_*$ of uniform variables attached to one member can be bounded in terms of $r,\eta$ alone. The integer $d$ enumerates the $P$ constructions and is never sampled. In particular, the second expression in~\eqref{eq:rounding-marginals} is the arithmetic mean of $P$ separately defined probabilities.
\end{proposition}

We prove \textbf{Proposition~\ref{prop:rounding}} in three steps.
The first step assumes a positive lower bound on the weights and a finite upper bound on a positive coefficient $z_a$ assigned to each member.
In \textbf{Lemma~\ref{lem:bounded-rounding}}, we bound the expected value of the weighted sum $\sum_{a\in J}z_a$ for the resulting matching $J$.
The second step removes the lower bound on weights by Poisson sampling.
The third step chooses the positive coefficients $z_a$ so that the average lower bound in~\eqref{eq:rounding-marginals} holds for every member separately. Once chosen from the given input, the coefficients $z_a$ are kept fixed in all subsequent experiments.

\subsection{Weights bounded below and coefficients bounded above}

\begin{lemma}\label{lem:bounded-rounding}
Fix integers $r\ge2$ and $M\ge1$, and reals $R\ge1$ and $0<\sigma<1$.
Let $\mathcal A$ have $n$ members of size at most $r$.
Suppose that $\frac{1}{M}\le h_a\le1$, $\sum_{a\ni v}h_a\le1$ for every $v$, and $1\le z_a\le R$ for every $a$.
There are a positive integer $T$ and random matchings $J_s$, $T\le s<2T$, such that
\begin{equation}\label{eq:bounded-rounding-conclusion}
\mathbb P(a\in J_s)\le\frac{h_a}{C},
\qquad
\frac{1}{T}\sum_{s=T}^{2T-1}\mathbb E\sum_{a\in J_s}z_a
\ge\frac{1-\sigma}{C}\sum_{a\in\mathcal A}z_ah_a.
\end{equation}

There are also positive integers $R_0,q_0$, depending only on $r,M,R,\sigma$, such that every member receives $q_0$ independent uniform variables. For each $a$, the indicator $\mathbf1_{\{a\in J_s\}}$ uses only $B_{U(\mathcal A)}(a,R_0)$ and the variables attached to members of $B_{U(\mathcal A)}(a,R_0)$. More precisely, two collections of uniform values agreeing on $B_{U(\mathcal A)}(a,R_0)$ give the same value of $\mathbf1_{\{a\in J_s\}}$. The integers $T,R_0,q_0$ and all selection functions are chosen before the uniform variables are drawn.
\end{lemma}

\begin{proof}
If $n=0$, take the empty matching.
Assume that $n>0$.
Since every positive weight is at least $\frac{1}{M}$, at most $M$ members contain a specified vertex.
Thus $U(\mathcal A)$ has maximum degree at most $B=rM$.
All constants in the proof of \textbf{Lemma~\ref{lem:bounded-rounding}} are chosen in terms of $r,M,R,\sigma$.

\par\smallskip\emph{The matrix and the iterates.}
Fix an integer $L\ge1$.
For each $a\in\mathcal A$ and $j\in[L]$, introduce a column $(a,j)$.
Let $\mathcal C=\mathcal A\times[L]$ be the set of column labels, so there are $N=nL$ columns.
The columns $(a,j)$ with a fixed value of $j$ form layer $j$.
For each vertex $v$ appearing in the family and each $j\in[L]$, introduce a row $(v,j)$, containing a $1$ in column $(a,j)$ exactly when $v\in a$; its right-hand side is $1$.
For each member $a$, introduce one further row containing a $1$ in every column $(a,j)$, $j\in[L]$; its right-hand side is $Lh_a$.
The vertex rows and member rows have distinct labels. Let $\mathcal R$ be the set of all row labels. All other entries are zero. Write $A=(A_{i,c})_{i\in\mathcal R,c\in\mathcal C}$ for the matrix of the vertex and member rows, and let $b=(b_i)_{i\in\mathcal R}$ be the vector of right-hand sides. A nonnegative column vector $y$ is called feasible if $Ay\le b$ coordinatewise.
The vector $x^0$ defined by $x^0_{a,j}=h_a$ satisfies $Ax^0\le b$.
There are at most $(r+1)N$ rows.
Every column contains at most $r+1$ ones, and every row contains at most $D_*:=\max\{M,L\}$ ones.

Fix $0<\mu<1$ and independently sample $U_{a,j}$ uniformly from $[0,1]$.
Put $w_{a,j}=z_a(1+\mu U_{a,j})$ and $W=2R$, so $1\le w_{a,j}\le W$.
For a parameter $\Lambda\ge1$, put $H=\Lambda(r+1)D_*$ and define
\begin{align}
F(u)&=w\cdot u-\frac{\Lambda}{2}\sum_i((Au)_i-b_i)_+^2,
\label{eq:penalized-objective}\\
p_i(u)&=\Lambda((Au)_i-b_i)_+,\qquad u^0=0,\notag\\
u^{s+1}&=\Pi_{[0,2]^N}\left(u^s+\frac{w-A^{\mathsf T}p(u^s)}{H}\right)\qquad(s\ge0).
\label{eq:gradient-iteration}
\end{align}

Here $(q)_+=\max\{q,0\}$ and the Euclidean projection onto $[0,2]^N$ satisfies $(\Pi_{[0,2]^N}(z))_c=\min\{2,\max\{0,z_c\}\}$ for every coordinate $c$.
For every real column vector $q$, Cauchy--Schwarz and the row and column bounds give
\begin{equation}\label{eq:matrix-norm-bound}
\|Aq\|_2^2
=\sum_i\left(\sum_c A_{ic}q_c\right)^2
\le D_*\sum_i\sum_c A_{ic}q_c^2
\le(r+1)D_*\|q\|_2^2.
\end{equation}

Thus $\|A\|_2^2\le(r+1)D_*$.
The real function $t\mapsto(t-b_i)_+^2$ has the nondecreasing derivative $2(t-b_i)_+$ and is therefore convex. The composition $u\mapsto((Au)_i-b_i)_+^2$ is also convex. Thus the function $F$ in~\eqref{eq:penalized-objective} is concave.
For all real $q,q'$, we have $|(q-b_i)_+-(q'-b_i)_+|\le|q-q'|$.
Using~\eqref{eq:matrix-norm-bound} in the formula for $\nabla F$ therefore gives $\|\nabla F(u)-\nabla F(v)\|_2\le\Lambda\|A\|_2^2\|u-v\|_2\le H\|u-v\|_2$.

For one update $u^+$ in~\eqref{eq:gradient-iteration}, \textbf{Lemma~\ref{lem:projection}} gives
\begin{equation}\label{eq:update-inner-product}
(\nabla F(u)-H(u^+-u))\cdot(v-u^+)\le0
\qquad(v\in[0,2]^N).
\end{equation}

Combining~\eqref{eq:update-inner-product} with concavity and $\|\nabla F(u)-\nabla F(v)\|_2\le H\|u-v\|_2$ gives
\begin{align}
F(u^+)-F(u)&\ge\frac{H}{2}\|u^+-u\|_2^2,\label{eq:gradient-increase}\\
F(v)-F(u^+)&\le\frac{H}{2}\bigl(\|v-u\|_2^2-\|v-u^+\|_2^2\bigr).
\label{eq:gradient-telescope}
\end{align}

Indeed, integrating the gradient along the segment from $u$ to $u^+$ and using $\|\nabla F(u)-\nabla F(v)\|_2\le H\|u-v\|_2$ gives $F(u^+)\ge F(u)+\nabla F(u)\cdot(u^+-u)-\frac{H}{2}\|u^+-u\|_2^2$.
For~\eqref{eq:gradient-increase}, use~\eqref{eq:update-inner-product} with $v=u$.
For~\eqref{eq:gradient-telescope}, combine $F(v)\le F(u)+\nabla F(u)\cdot(v-u)$ with $F(u^+)\ge F(u)+\nabla F(u)\cdot(u^+-u)-\frac{H}{2}\|u^+-u\|_2^2$, and then apply~\eqref{eq:update-inner-product} with $v$.

Apply~\eqref{eq:gradient-telescope} with $v=x^0$ to the updates from $u^0$ through $u^T$ and sum. The squared-distance terms telescope. Since~\eqref{eq:gradient-increase} makes $F(u^s)$ nondecreasing, each of $F(u^1),\ldots,F(u^T)$ is at most $F(u^T)$. Hence $T(F(x^0)-F(u^T))\le\frac{H}{2}\|x^0\|_2^2\le\frac{HN}{2}$. Using $F(x^0)=w\cdot x^0$ and $F(u^0)=0$, we obtain
\begin{align}
&F(u^s)\ge w\cdot x^0-\frac{HN}{2T},\qquad F(u^s)\ge0\quad(s\ge T),\label{eq:iterate-objective}\\
&\frac{1}{T}\sum_{s=T}^{2T-1}\|u^{s+1}-u^s\|_2^2\le\frac{4WN}{HT}.
\label{eq:iterate-changes}
\end{align}

To prove~\eqref{eq:iterate-changes}, sum~\eqref{eq:gradient-increase} over $s=T,\ldots,2T-1$. The objective differences telescope, giving $\frac{H}{2}\sum_{s=T}^{2T-1}\|u^{s+1}-u^s\|_2^2\le F(u^{2T})-F(u^T)\le2WN$, because $F(u^T)\ge0$ and $F(u)\le2WN$ on $[0,2]^N$. Dividing by $\frac{HT}{2}$ proves~\eqref{eq:iterate-changes}. For each $s\in\{T,\ldots,2T-1\}$, we construct $J_s$ from $u^s$, the vector after the $s$th update.

\par\smallskip\emph{Satisfying the matrix constraints.}
Fix $s\in\{T,\ldots,2T-1\}$ and write $u=u^s$, $p=p(u)$, and $e_i=((Au)_i-b_i)_+$.
By $F(u)\ge0$ and Cauchy--Schwarz,
\begin{equation}\label{eq:constraint-excess}
\sum_i e_i^2\le\frac{4WN}{\Lambda},
\qquad
\sum_i e_i\le N\sqrt{\frac{4W(r+1)}{\Lambda}}.
\end{equation}

For each column $c$, define
\begin{equation}\label{eq:feasible-scaling}
\lambda_c=\max\left\{1,\max_{i:A_{ic}=1}\frac{(Au)_i}{b_i}\right\},
\qquad
\widetilde y_c=\frac{u_c}{\lambda_c}.
\end{equation}

Equation~\eqref{eq:feasible-scaling} gives $A\widetilde y\le b$. Indeed, if $(Au)_i\le b_i$, then the $i$th row inequality follows from $0\le\widetilde y\le u$. If $(Au)_i>b_i$, equation~\eqref{eq:feasible-scaling} gives $\lambda_c\ge\frac{(Au)_i}{b_i}$ for each column with $A_{ic}=1$, so $\sum_c A_{ic}\widetilde y_c\le\frac{b_i}{(Au)_i}\sum_c A_{ic}u_c=b_i$.
Moreover, $\sum_c(u_c-\widetilde y_c)\le\sum_i e_i$.
To prove $\sum_c(u_c-\widetilde y_c)\le\sum_i e_i$, assign each decreased coordinate $c$ to a row $i$ attaining the maximum in~\eqref{eq:feasible-scaling}, breaking ties by a fixed order.
For a fixed row $i$ with positive excess, scaling all its coordinates by $\frac{b_i}{(Au)_i}$ decreases their total by exactly $e_i$.
The decrease of the coordinates assigned to that row is at most $e_i$.
Summing over rows gives $\sum_c(u_c-\widetilde y_c)\le\sum_i e_i$.

Choose $0<\beta<\frac{1}{2}$ and $0<\tau<1$.
Set $y_c=0$ if $u_c<\beta$, if $u_c>\frac{3}{2}$, or if $|w_c-(A^{\mathsf T}p)_c|>\tau$; otherwise set $y_c=\widetilde y_c$.
Since $0\le y\le\widetilde y$ and $A\widetilde y\le b$ by~\eqref{eq:feasible-scaling}, we have $Ay\le b$.
Every positive coordinate of $y$ satisfies
\begin{equation}\label{eq:small-coordinate-gradient}
|w_c-(A^{\mathsf T}p)_c|\le\tau.
\end{equation}

In~\eqref{eq:small-coordinate-gradient}, the vector $p=p(u)$ is calculated before the scaling in~\eqref{eq:feasible-scaling}.

There are three additional losses in passing from $\widetilde y$ to $y$.
The coordinates with $u_c<\beta$ have total at most $\beta N$.
If $u_c>\frac{3}{2}$, any vertex row containing $c$ has excess greater than $\frac{1}{2}$.
Assign each such column to one incident vertex row.
Each vertex row receives at most $M$ columns, so at most $2M\sum_i e_i$ columns are assigned in this way; their coordinates sum to at most $4M\sum_i e_i$.
Finally, if $\beta\le u_c\le\frac{3}{2}$ and~\eqref{eq:small-coordinate-gradient} fails, the next update changes coordinate $c$ by at least $\gamma:=\min\{\beta,\frac{1}{2},\frac{\tau}{H}\}$.
By~\eqref{eq:iterate-changes}, the average number of coordinates satisfying $\beta\le u_c\le\frac{3}{2}$ and violating~\eqref{eq:small-coordinate-gradient}, over $s=T,\ldots,2T-1$, is at most $\frac{4WN}{HT\gamma^2}$.
Each coordinate is at most $2$.
Combining~\eqref{eq:constraint-excess} with the losses from $u_c<\beta$, from $u_c>\frac{3}{2}$, and from violations of~\eqref{eq:small-coordinate-gradient}, we obtain
\begin{equation}\label{eq:feasibility-loss}
\frac{1}{T}\sum_{s=T}^{2T-1}\sum_c(u_c^s-y_c^s)\le NE_1,
\qquad
E_1=(1+4M)\sqrt{\frac{4W(r+1)}{\Lambda}}+\beta+\frac{8W}{HT\gamma^2}.
\end{equation}

Each vertex row has right-hand side $1$, and for every column $c$ there is a vertex row $i$ with $A_{ic}=1$.
Thus $0\le y_c\le1$.

\par\smallskip\emph{A system of distinct representatives for the column supports.}
Join two columns when they have a common row.
The graph on columns joined by common rows has maximum degree at most $D_c=(r+1)D_*$.
Fix $0<\xi_0<\frac{1}{2}$.
Let $S$ be a fixed set of $d$ columns whose incident rows number $l<(1-\xi_0)d$.
If all columns in $S$ satisfy~\eqref{eq:small-coordinate-gradient}, then $w_S$ lies within Euclidean distance $\tau\sqrt d$ of the fixed subspace
$V_S=\{(A^{\mathsf T}q)_S:q\in\mathbb R^{\mathcal R}\}$, whose dimension is at most $l$.
Although $p$ depends on all perturbations, the requirement that $w_S$ lie within Euclidean distance $\tau\sqrt d$ of $V_S$, necessary for~\eqref{eq:small-coordinate-gradient}, involves only $w_S$ and the fixed matrix $A$.

Since $\|w_S\|_2\le W\sqrt d$ and $\tau<1$, every point of $V_S$ at distance at most $\tau\sqrt d$ from $w_S$ lies in the ball of radius $(W+1)\sqrt d$ about zero in $V_S$.
The ball of radius $(W+1)\sqrt d$ in $V_S$ has a $\tau\sqrt d$-net with at most $\left(1+\frac{2(W+1)}{\tau}\right)^l$ points.
Here a net means a set such that every point of the ball is at distance at most $\tau\sqrt d$ from a point in the net.
To obtain the cardinality bound, choose a maximal set of points in the radius-$(W+1)\sqrt d$ ball whose pairwise distances exceed $\tau\sqrt d$. Maximality makes the chosen set a $\tau\sqrt d$-net. The balls of radius $\frac{\tau\sqrt d}{2}$ about the chosen points are disjoint and lie in the ball of radius $(W+1+\frac{\tau}{2})\sqrt d$ in $V_S$. Comparing volumes in $V_S$, whose dimension is at most $l$, gives the bound $\left(1+\frac{2(W+1)}{\tau}\right)^l$.
A Euclidean ball of radius $2\tau\sqrt d$ in $\mathbb R^d$ has volume at most $(c_0\tau)^d$, where one may take $c_0=\sqrt{8\pi\mathrm e}$. Indeed, if $B$ is such a ball centered at zero, then $\exp\left(-\frac{\|x\|_2^2}{8\tau^2}\right)\ge\exp\left(-\frac{d}{2}\right)$ for $x\in B$. Integrating and using the product of the one-dimensional Gaussian integrals gives
\[
\operatorname{vol}(B)
\le \mathrm e^{\frac{d}{2}}\int_{\mathbb R^d}
       \exp\left(-\frac{\|x\|_2^2}{8\tau^2}\right)\,dx
=\mathrm e^{\frac{d}{2}}(\sqrt{8\pi}\,\tau)^d
=(c_0\tau)^d.
\]

The coordinates of $w_S$ are independent and each has density at most $\frac{1}{\mu}$, because $z_a\ge1$.
The union bound over the net therefore gives
\begin{equation}\label{eq:perturbation-deficiency-probability}
\mathbb P\bigl(\text{every }c\in S\text{ satisfies }\eqref{eq:small-coordinate-gradient}\bigr)
\le\left[\frac{(1+2(W+1))c_0\tau^{\xi_0}}{\mu}\right]^d.
\end{equation}

To derive~\eqref{eq:perturbation-deficiency-probability}, we used $0<\tau<1$ and $d-l>\xi_0d$ to bound $\left(1+\frac{2(W+1)}{\tau}\right)^l(c_0\tau)^d\mu^{-d}$ by $\left[\frac{(1+2(W+1))c_0\tau^{\xi_0}}{\mu}\right]^d$.

The number of connected sets of $d$ columns containing a specified column is at most $(1+D_c)^{2d}$.
Indeed, a spanning tree on such a set has a closed traversal of length $2(d-1)$ starting at that column.
Fix $\omega>0$.
After $L,\mu,\xi_0,R,M$ have been fixed, choose $\tau$ so small that $\sum_{d\ge1}(1+D_c)^{2d}\left[\frac{(1+2(W+1))c_0\tau^{\xi_0}}{\mu}\right]^d\le\omega$. The summands combine the bound on connected sets with~\eqref{eq:perturbation-deficiency-probability}.

For the analysis only, mark a column $c$ if there is a connected set $S$ of columns containing $c$ such that every member of $S$ satisfies~\eqref{eq:small-coordinate-gradient} and $|\{i:A_{i,b}=1\text{ for some }b\in S\}|<(1-\xi_0)|S|$.
If $Z$ is the number of marked columns, the union bound and~\eqref{eq:perturbation-deficiency-probability}, with the chosen $\tau$, give $\mathbb E Z\le\omega N$ at each fixed time.
No step of the random construction needs to identify marked columns.
Consider any set $S$ of unmarked columns $c$ with $y_c>0$.
The connected components of $S$ meet disjoint sets of rows.
Every component $S'$ has at least $(1-\xi_0)|S'|$ incident rows: otherwise every column of $S'$ would be marked, since the positive coordinates of $y$ satisfy~\eqref{eq:small-coordinate-gradient}.
Summing over the components, $S$ has at least $(1-\xi_0)|S|$ incident rows.

Apply \textbf{Corollary~\ref{cor:hall-deficiency}} with $b=\xi_0N$ to the bipartite graph whose left vertices are the unmarked columns $c$ with $y_c>0$, whose right vertices are the rows, and whose edges are the pairs $ci$ with $A_{ic}=1$.
The resulting matching gives a system of distinct representatives for the supports of all but at most $\xi_0N$ unmarked columns $c$ with $y_c>0$.
At most $n$ columns are matched to the rows indexed by members of $\mathcal A$.
Discard the columns assigned to member rows as well.
All remaining columns are matched to distinct vertex rows $(v,j)$.
Consequently, for each layer $j$ there is an integer $b_j$ such that deleting $b_j$ columns with $y_{a,j}>0$ leaves a family of sets $a$ having a system of distinct representatives, and
\begin{equation}\label{eq:layer-deletion-bound}
\sum_{j=1}^L b_j\le Z+\xi_0N+n.
\end{equation}

The number of unmatched columns is an integer at most $\xi_0N$. Thus~\eqref{eq:layer-deletion-bound} needs no additional rounding term.

\par\smallskip\emph{Removing members in successive rounds.}
Fix a layer $j$ and write $y_a=y_{a,j}$.
For a subfamily $\mathcal F\subseteq\{a\in\mathcal A:y_{a,j}>0\}$, let
$W_{\mathcal F}(a)=\sum_{b\in\mathcal F:\,b\cap a\ne\varnothing}y_b$, including $b=a$.
Deleting at most $b_j$ members from $\mathcal F$ leaves a family having a system of distinct representatives. Thus there is an injective map from a subfamily of size at least $|\mathcal F|-b_j$ into $\bigcup\mathcal F$, assigning to each member a vertex contained in that member. Hence $|\bigcup\mathcal F|\ge|\mathcal F|-b_j$.
Let $l$ be the number of vertices in $\bigcup\mathcal F$ that belong to at most $r$ members of $\mathcal F$.
Every other vertex belongs to at least $r+1$ members, and every member has size at most $r$.
Counting incidences between $\mathcal F$ and $\bigcup\mathcal F$ gives
\begin{equation}\label{eq:few-incident-members}
r|\mathcal F|\ge\sum_v d_{\mathcal F}(v)
\ge(r+1)|\bigcup\mathcal F|-rl,
\qquad
l\ge\frac{|\mathcal F|-(r+1)b_j}{r}.
\end{equation}

At each vertex belonging to at most $r$ members of $\mathcal F$, choose an incident member $a$ of maximum weight.
Suppose first that $y_a\ge\frac{1}{r}$. Summing the vertex constraints over the vertices of $a$ counts every intersecting member at least once and counts $a$ exactly $|a|$ times. Subtracting the $|a|-1$ extra contributions from $a$ therefore gives $W_{\mathcal F}(a)\le |a|-(|a|-1)y_a\le |a|-\frac{|a|-1}{r}\le r-1+\frac{1}{r}=C$, using $y_a\ge\frac{1}{r}$ and $|a|\le r$.
If $y_a<\frac{1}{r}$, the sum at the chosen vertex is at most $ry_a$, and at every other vertex of $a$ the sum excluding $a$ is at most $1-y_a$.
Hence $W_{\mathcal F}(a)\le ry_a+(|a|-1)(1-y_a)\le r-1+y_a\le C$.
A member can be chosen at at most $r$ vertices.
By~\eqref{eq:few-incident-members}, at least $\frac{|\mathcal F|-(r+1)b_j}{r^2}$ members therefore satisfy $W_{\mathcal F}(a)\le C$.

In every layer perform the following operations, using the same integers $J,Q$.
\begin{enumerate}[label=\textup{(\arabic*)}]
\item Begin with the members $a$ for which $y_{a,j}>0$.
For $J$ rounds, simultaneously remove every member $a$ for which $W_{\mathcal F}(a)\le C$ in the family present at the start of the round.
Record the removal round of each removed member $a$.
Discard members still present after round $J$.
\item Give each removed member an independent uniform label in $[Q]$.
Discard a removed member $a$ if $a$ intersects another member with the same removal round and the same label.
\end{enumerate}

Let $n_j^{(h)}$ be the number of members remaining in layer $j$ after $h$ rounds. By~\eqref{eq:few-incident-members} and the fact that a member is chosen at at most $r$ vertices, at least $\frac{n_j^{(h)}-(r+1)b_j}{r^2}$ members are removed in round $h+1$. Hence $n_j^{(h+1)}\le\left(1-\frac{1}{r^2}\right)n_j^{(h)}+\frac{(r+1)b_j}{r^2}$ for $0\le h<J$. Iterating the recurrence from $n_j^{(0)}\le n$ and summing the geometric progression gives $n_j^{(J)}\le\left(1-\frac{1}{r^2}\right)^Jn+(r+1)b_j$.
For a removed member $a$, the probability of meeting another member with the same removal round and label is at most $\frac{B}{Q}$, by the union bound over the neighbors of $a$ in $U(\mathcal A)$.
Members having the same removal round and label are pairwise disjoint after the discards.
Taking expectations in~\eqref{eq:layer-deletion-bound}, the expected number of discarded columns $c$ with $y_c>0$, divided by $N$, is at most
\begin{equation}\label{eq:removal-loss}
E_2=\left(1-\frac{1}{r^2}\right)^J+(r+1)\left(\omega+\xi_0+\frac{1}{L}\right)+\frac{B}{Q}.
\end{equation}

\par\smallskip\emph{Choosing a matching.}
Let $\mathscr P$ be the random vector consisting of the column perturbations and the random labels in all layers. The independent uniform variables used for inclusion decisions are not coordinates of $\mathscr P$. Condition on $\mathscr P$, so the iterates, the vectors $y$, and the removal rounds are fixed.
In each layer process the surviving members in decreasing order of removal round and, within a round, in increasing order of label.
The members of each resulting group are pairwise disjoint.
Put $C'=C+\mu$.
We define decisions for successive groups using fresh independent uniform variables, one for each member in each layer.
When member $a$ is processed, let $\pi_a$ be the probability, conditional on $\mathscr P$ alone, that no previously included member intersects $a$.
The conditional probability $\pi_a$ averages over the fresh inclusion variables of all groups processed before $a$, without conditioning on the realized values of those variables.

Every previously processed member intersecting $a$ belonged to the family from which $a$ was removed.
Suppose inductively that each such member $b$ has inclusion probability $\frac{y_b}{C'}$ conditional on $\mathscr P$.
Since $W_{\mathcal F}(a)\le C$ when $a$ is removed, the union bound gives
\begin{equation}\label{eq:free-probability}
\pi_a\ge1-\frac{C-y_a}{C'}\ge\frac{y_a}{C'}>0.
\end{equation}

If no previously included member intersects $a$ in the realized partial matching, include $a$ with probability $\frac{y_a}{C'\pi_a}$, using the fresh inclusion variable attached to $a$.
Otherwise do not include $a$.
Equation~\eqref{eq:free-probability} gives $\frac{y_a}{C'\pi_a}\in[0,1]$.
The inclusion probability of $a$, conditional on $\mathscr P$, is exactly $\pi_a\frac{y_a}{C'\pi_a}=\frac{y_a}{C'}$, completing the induction.
Since members with the same removal round and label are pairwise disjoint, the decisions in each group may be made simultaneously.

For each fixed time index, the construction of $\mathscr P$, the calculations of $y$, the removal operations, and the inclusion rules in~\eqref{eq:free-probability} are invariant under permutations of the $L$ layers.
Let $J_s$ be the matching obtained in layer $1$ at time index $s$.
The row indexed by $a$ gives $\sum_j y_{a,j}\le Lh_a$ for every realization.
The conditional marginal formula and layer symmetry therefore give $\mathbb P(a\in J_s)\le \frac{h_a}{C'}\le \frac{h_a}{C}$.

Put $S_h=\sum_a z_ah_a$, so $S_h\ge \frac{n}{M}$.
By~\eqref{eq:iterate-objective}, $\sum_c w_cu_c\ge LS_h-\frac{HN}{2T}$ at each time under consideration.
The average loss when replacing $u$ by $y$ is at most $WNE_1$, by~\eqref{eq:feasibility-loss}.
The additional expected loss from discarding columns $c$ with $y_c>0$ is at most $WNE_2$, by~\eqref{eq:removal-loss} and $y_c\le1$.
Finally, $z_a\ge \frac{w_{a,j}}{1+\mu}$, and each surviving column is included with conditional probability $\frac{y_{a,j}}{C'}$.
Average over the layers and time indices to obtain
\begin{equation}\label{eq:bounded-weighted-sum-estimate}
\frac{1}{T}\sum_{s=T}^{2T-1}\mathbb E\sum_{a\in J_s}z_a
\ge\frac{S_h-n[\frac{H}{2T}+W(E_1+E_2)]}{C'(1+\mu)}.
\end{equation}

Choose the parameters in the following order.
\begin{enumerate}[label=\textup{(\roman*)}]
\item Choose $0<\mu<1$ such that $\frac{C}{C'(1+\mu)}>1-\frac{\sigma}{3}$. The choice is possible because $C'=C+\mu$ and $\frac{C}{C'(1+\mu)}\to1$ as $\mu\to0$.
\item Choose $\omega>0$, $0<\xi_0<\frac{1}{2}$, and positive integers $L,J,Q$ such that $WE_2\le\frac{\sigma}{12M}$, using the formula for $E_2$ in~\eqref{eq:removal-loss}. Each term in that formula can be made arbitrarily small by decreasing $\omega,\xi_0$ or increasing $L,J,Q$, respectively.
\item Choose $0<\beta<\frac{1}{2}$ and $\Lambda\ge1$ such that $W[(1+4M)\sqrt{\frac{4W(r+1)}{\Lambda}}+\beta]\le\frac{\sigma}{12M}$.
\item Choose $\tau$ such that $\sum_{d\ge1}(1+D_c)^{2d}\left[\frac{(1+2(W+1))c_0\tau^{\xi_0}}{\mu}\right]^d\le\omega$, as required when applying~\eqref{eq:perturbation-deficiency-probability}. With $H$ and $\gamma$ now fixed, choose an integer $T$ satisfying $T\ge\frac{6M}{\sigma}\left[\frac{H}{2}+\frac{8W^2}{H\gamma^2}\right]$.
\end{enumerate}

Substituting the chosen $\Lambda,\beta,T$ into~\eqref{eq:feasibility-loss} and using $WE_2\le\frac{\sigma}{12M}$ from~\eqref{eq:removal-loss}, we obtain $\frac{H}{2T}+W(E_1+E_2)\le\frac{\sigma}{3M}$.
Since $n\le MS_h$, equation~\eqref{eq:bounded-weighted-sum-estimate} is at least $\frac{(1-\frac{\sigma}{3})^2S_h}{C}\ge\frac{(1-\sigma)S_h}{C}$.
Combining~\eqref{eq:bounded-weighted-sum-estimate} with $\mathbb P(a\in J_s)\le\frac{h_a}{C}$, obtained from the conditional marginal formula and $\sum_j y_{a,j}\le Lh_a$, proves~\eqref{eq:bounded-rounding-conclusion}.

We now verify that $\mathbf1_{\{a\in J_s\}}$ is determined by the random variables attached to $B_{U(\mathcal A)}(a,R_0)$, as required in \textbf{Lemma~\ref{lem:bounded-rounding}}. A single update in~\eqref{eq:gradient-iteration} uses only columns sharing a row. Columns sharing a row belong either to the same member or to two intersecting members. After $s$ updates, the coordinates belonging to $a$ therefore use perturbations at members of $B_{U(\mathcal A)}(a,s)$. Scaling in~\eqref{eq:feasible-scaling} and the test~\eqref{eq:small-coordinate-gradient} use the variables at members of one additional neighborhood. Since $s<2T$, each value $y_{a,j}$ is determined within radius $2T$.

The $J$ removal rounds use at most $J$ additional neighborhoods, and the comparison of labels uses one more. Thus all perturbations, labels, and quantities computed from them that are needed for the decision at $a$ are determined within radius $2T+J+1$. The marked columns and the systems of distinct representatives are used only in the proof of~\eqref{eq:layer-deletion-bound}; neither is used to construct $J_s$.

There are at most $JQ$ groups in the final inclusion procedure. Conditional on $\mathscr P$, decisions in the first group use only the inclusion variables of members in that group. Inductively, the decision for $a$ in the next group uses the earlier decisions at neighbors of $a$ and the inclusion variable attached to $a$. The probability $\pi_a$ in~\eqref{eq:free-probability} is obtained by integrating the inclusion variables used for those earlier decisions, so $\pi_a$ is determined by the coordinates of $\mathscr P$ needed for the earlier decisions. Processing each new group increases the required radius by at most one. Starting from radius $2T+J+1$ for the perturbations and labels and allowing $JQ$ additional neighborhoods, we may take $R_0=2T+J+JQ+3$. Finally, the perturbation, label, and inclusion choices require one independent uniform variable per member in each layer, so $q_0=3L$ suffices. The integers $T,R_0,q_0$ depend only on $r,M,R,\sigma$, as required.
\end{proof}

\subsection{Removing the lower bound on positive weights}

\begin{lemma}\label{lem:poisson-expectation}
Let $0<\delta<\frac{1}{4}$ and suppose that $\sum_{a\ni v}x_a\le1$ for every $v$, with $|a|\le r$.
For sufficiently large $m=m(r,\delta)$, put $M=\lceil(1+\delta)m\rceil$, sample the variables of \textbf{Definition~\ref{def:local-random-matching}}, and set $h_a=\frac{K_a}{M}$ for $a\in\mathcal A_K$ and $h_a=0$ for $a\notin\mathcal A_K$.
Then
\begin{equation}\label{eq:poisson-weight-expectation}
\sum_{a\ni v}h_a\le1\quad\text{for every }v,
\qquad
(1-3\delta)x_a\le\mathbb Eh_a\le x_a\quad\text{for every }a.
\end{equation}

Each positive $h_a$ lies in $[\frac{1}{M},1]$.
\end{lemma}

\begin{proof}
If any member incident with $v$ is retained, condition~\eqref{eq:retained-family} gives $T_v\le M$.
Consequently $\sum_{a\in\mathcal A_K:a\ni v}K_a\le M$.
If no member containing $v$ is retained, then $\sum_{a\ni v}h_a=0$.
In both cases, $\sum_{a\ni v}h_a\le1$, proving the first inequality in~\eqref{eq:poisson-weight-expectation}.
The upper expectation bound follows from $h_a\le \frac{K_a}{M}$ and $\mathbb E\frac{K_a}{M}=\frac{mx_a}{M}\le x_a$.

For every bounded function $f:\mathbb Z_{\ge0}^{\mathcal A}\to\mathbb R$, the Poisson mass function gives $\mathbb E[K_af(K)]=mx_a\mathbb E[f(K+\mathbf e_a)]$, where $\mathbf e_a$ has coordinate $1$ at $a$ and coordinate $0$ at every $b\ne a$.
Apply $\mathbb E[K_af(K)]=mx_a\mathbb E[f(K+\mathbf e_a)]$ with $f(K)=\mathbf1_{\{T_v\le M\text{ for every }v\in a\}}$.
For the modified vector $K+\mathbf e_a$, coordinate $a$ is positive automatically.
Before adding $\mathbf e_a$, each sum $T_v$ with $v\in a$ is Poisson with mean at most $m$.
Therefore the union bound over at most $r$ vertices gives
\begin{equation}\label{eq:poisson-loss-bound}
\mathbb Eh_a\ge\frac{mx_a}{M}\left[1-r\mathbb P(\operatorname{Pois}(m)>M-1)\right].
\end{equation}

For fixed $\delta>0$, the probability in~\eqref{eq:poisson-loss-bound} tends to zero exponentially as $m\to\infty$, while $\frac{m}{M}\to\frac{1}{1+\delta}$.
For example, the exponential Markov inequality applied to a Poisson variable $X$ gives $\mathbb P(X\ge q)\le\exp(-q\log(\frac{q}{m})+q-m)$ for $q>m$.
Choose $m$ so large that $\frac{m}{M}[1-r\mathbb P(\operatorname{Pois}(m)>M-1)]\ge1-3\delta$.
Substituting into~\eqref{eq:poisson-loss-bound} proves the lower expectation bound in~\eqref{eq:poisson-weight-expectation}.
The bound on positive weights follows directly from retention and $K_a\ge1$.
\end{proof}

\begin{proof}[Proof of \textbf{Proposition~\ref{prop:rounding}}]
For $\mathcal A=\emptyset$, take $I_d=\emptyset$ for every $d$; we may therefore assume that $n_0=|\mathcal A|>0$.
Choose $0<\delta<\frac{1}{4}$ such that $\kappa_\delta\ge1-\frac{\eta}{2}$, where $\kappa_\delta=\frac{(1-\delta)^3(1-3\delta)}{1+\delta}$, and put $\zeta=\frac{\eta}{2}$. This choice is possible because $\kappa_\delta\to1$ as $\delta\to0$. Choose $m,M$ by \textbf{Lemma~\ref{lem:poisson-expectation}}.
We first give a construction for every fixed vector of positive coefficients $z=(z_a:a\in\mathcal A)$. The constants in the construction and the integer $P$ will be independent of the values of $z_a$.

Independently sample a uniform variable $V_a\in[0,1]$ for each member and put $\widehat z_a=z_a(1+\delta V_a)$.
The subsequent construction uses $\widehat z$, and does not otherwise use $z$.
Assign integer $t_a=\lfloor\log_2\widehat z_a\rfloor$ to member $a$.
Choose positive integers $g,\ell$ and put $q=g+\ell$.
For each fixed integer $c\in\{0,\ldots,q-1\}$, delete members whose integer $t_a$ is congruent to one of $c,c+1,\ldots,c+g-1$ modulo $q$.
The remaining integers form maximal intervals of $\ell$ consecutive integers, separated by $g$ deleted integers; we call each such interval a \emph{block}.
In each nonempty block, divide the coefficients $\widehat z_a$ of the members in the block by $2^{t_0}$, where $t_0$ is the least integer of the full block of $\ell$ consecutive integers determined by the fixed shift $c$, whether or not a member is assigned to $t_0$.
The coefficients in each block, after division by $2^{t_0}$, lie in $[1,2^\ell]$.
Apply \textbf{Lemma~\ref{lem:bounded-rounding}} to the retained members of that block with weights $h_a$, coefficient bound $R=2^\ell$, and error $\sigma=\delta$.
Use the same parameters and time indices in every block.

For a fixed shift $c$ and time index $s$, let $\mathcal J$ be the union of the matchings obtained from \textbf{Lemma~\ref{lem:bounded-rounding}} in the nonempty blocks. From $\mathcal J$, discard $a$ whenever $a$ intersects a member of $\mathcal J$ from a block with larger integers.
When deciding whether to discard $a$, compare $a$ with all members of $\mathcal J$ from blocks with larger integers, including members that are themselves discarded.
Let $I_{c,s}$ be the surviving family. Then $I_{c,s}$ is a matching: members of $\mathcal J$ in the same block are already disjoint, and if two members from different blocks intersect, the member in the lower block is discarded.
For each $a\in\mathcal J$, at most $B=rM$ members $b\in\mathcal J$ from lower blocks intersect $a$, since $\mathcal J\subseteq\mathcal A_K$. Each such $b$ satisfies $\widehat z_b\le2^{-g}\widehat z_a$.
Assign each discarded $b\in\mathcal J$ to one intersecting $a\in\mathcal J$ from a block with larger integers. Each $a$ receives at most $B$ assignments, each of weight at most $2^{-g}\widehat z_a$. Thus the sum of $\widehat z_b$ over discarded members is at most $B2^{-g}\sum_{a\in\mathcal J}\widehat z_a$.

Choose $g$ so that $B2^{-g}\le\delta$, and then choose $\ell$ so that $\frac{g}{g+\ell}\le\delta$.
Average uniformly over all $q$ shifts and all $T$ time indices supplied by \textbf{Lemma~\ref{lem:bounded-rounding}}.
There are $P=qT$ pairs $(c,s)$, each specifying one shift and one iteration index. Enumerate the pairs by $d\in[P]$, and write $I_d=I_{c,s}$ for the matching corresponding to pair $(c,s)$. Neither component of $(c,s)$ is sampled.
Each member is deleted for exactly $g$ shifts.
Conditionally on the Poisson variables and on $\widehat z$, equation~\eqref{eq:bounded-rounding-conclusion} gives $\mathbb P(a\in\mathcal J\mid K,\widehat z)\le\frac{h_a}{C}$. Since $I_{c,s}\subseteq\mathcal J$, the conditional probability that $a$ belongs to $I_{c,s}$ is also at most $\frac{h_a}{C}$.
The three factors $1-\delta$ in the averaged weighted-sum bound respectively account for the error in \textbf{Lemma~\ref{lem:bounded-rounding}}, the deleted shifts, and combining the blocks.
Put $p_{a,d}(z)=\mathbb P(a\in I_d)$ and $\overline p_a(z)=\frac{1}{P}\sum_{d=1}^Pp_{a,d}(z)$.
Since $z_a\le\widehat z_a\le(1+\delta)z_a$, taking expectations and using~\eqref{eq:poisson-weight-expectation} gives
\begin{equation}\label{eq:arbitrary-coefficient-bound}
p_{a,d}(z)\le\frac{x_a}{C},
\qquad
\sum_a z_a\overline p_a(z)\ge\frac{\kappa_\delta}{C}\sum_a z_ax_a.
\end{equation}

For each index, group all variables by the member labels $a$.
The group indexed by $a$ consists of $K_a$, $V_a$, and the column perturbations, layer labels, and fresh inclusion variables attached to $a$ in \textbf{Lemma~\ref{lem:bounded-rounding}}.
All variables in the groups indexed by $a\in\mathcal A$ can be drawn independently in advance, including the variables that the construction does not use in a given outcome.
Processing all blocks by the same rules preserves the radius $R_0$ from \textbf{Lemma~\ref{lem:bounded-rounding}}. To decide whether $a\in\mathcal J$ is discarded because $a$ meets a member of $\mathcal J$ from a block with larger integers, we determine membership in $\mathcal J$ for each neighbor of $a$ in $U(\mathcal A_K)$. Thus $R_*=R_0+1$ suffices after the Poisson retention step. The variable $V_a$ and the $3L$ uniform variables used in \textbf{Lemma~\ref{lem:bounded-rounding}} form a tuple of length $q_*=1+3L$ for each member. Discrete labels in $[Q]$ are obtained by dividing one uniform variable into $Q$ intervals of equal length. The input variables and selection functions therefore satisfy \textbf{Definition~\ref{def:local-random-matching}}, with $R_*,q_*$ depending only on $r,\delta$.

For the given finite input, the functions $\overline p_a(z)$ are continuous as functions of the positive coefficients $(z_a)_{a\in\mathcal A}$.
Here is a direct verification that also handles all integer comparisons in the construction.
The density of $\widehat z_a$ is $f_{z_a}(t)=\frac{1}{\delta z_a}\mathbf1_{[z_a,(1+\delta)z_a]}(t)$.
When $z_a'>0$ tends to $z_a>0$, the endpoints of $[z_a',(1+\delta)z_a']$ and the height $\frac{1}{\delta z_a'}$ converge to the corresponding values for $f_{z_a}$, so $\int_{\mathbb R}|f_{z_a}-f_{z_a'}|\,dt\to0$.
For fixed $\widehat z$, all remaining random variables have distributions independent of $z$ and all decision rules are fixed measurable functions of their inputs.
Consequently their conditional probability of including $a$, averaged over the indices, is one measurable function $g_a(\widehat z)\in[0,1]$, independent of $z$.
Integrating $g_a$ against the product densities for $z$ and $z'$ gives
\begin{equation}\label{eq:coefficient-continuity}
|\overline p_a(z)-\overline p_a(z')|
\le\frac{1}{2}\int_{\mathbb R^{n_0}}
\left|\prod_b f_{z_b}(t_b)-\prod_b f_{z_b'}(t_b)\right|\,dt
\le\frac{1}{2}\sum_b\int_{\mathbb R}|f_{z_b}-f_{z_b'}|\,dt.
\end{equation}

The first inequality in~\eqref{eq:coefficient-continuity} uses $0\le g_a\le1$ and the fact that both product densities integrate to $1$.
The second inequality in~\eqref{eq:coefficient-continuity} follows by telescoping the difference of the products and integrating the remaining density factors.
Since $\mathcal A$ is finite and each integral $\int_{\mathbb R}|f_{z_b}-f_{z_b'}|\,dt$ tends to zero, the right-hand side of~\eqref{eq:coefficient-continuity} tends to zero.

We now choose the coefficients as a function of the fixed input, before any experiment is drawn.
Let $\Sigma_{n_0}$ be the probability simplex defined in \textbf{Subsection~\ref{s:prelim}}.
For $u\in\Sigma_{n_0}$, set $z_a=\frac{u_a+\frac{\zeta}{n_0}}{x_a}$ and $v_a(u)=\frac{\overline p_a(z)}{x_a}$, using the positive value $\zeta=\frac{\eta}{2}$ chosen at the start of the proof.
Then $v$ is continuous and, by~\eqref{eq:arbitrary-coefficient-bound},
\begin{equation}\label{eq:simplex-coefficient-inequality}
0\le v_a(u)\le\frac{1}{C},
\qquad
u\cdot v(u)\ge\frac{\kappa_\delta}{C}-\frac{\zeta}{C}.
\end{equation}

Indeed, $\sum_a z_ax_a=1+\zeta$ and
$\sum_a z_a\overline p_a(z)=u\cdot v(u)+\frac{\zeta}{n_0}\sum_a v_a(u)$, where $\frac{\zeta}{n_0}\sum_a v_a(u)\le\frac{\zeta}{C}$ by~\eqref{eq:arbitrary-coefficient-bound}.

The map $u\mapsto\Pi_{\Sigma_{n_0}}(u-v(u))$ is a continuous self-map of the compact convex simplex.
\textbf{Theorem~\ref{thm:brouwer}} gives $u^*\in\Sigma_{n_0}$ satisfying $u^*=\Pi_{\Sigma_{n_0}}(u^*-v(u^*))$.
Equation~\eqref{eq:convex-minimizer-inequality} in \textbf{Lemma~\ref{lem:projection}} then gives
$v(u^*)\cdot(y-u^*)\ge0$ for every $y\in\Sigma_{n_0}$.
Taking $y=\mathbf e_a$ and using~\eqref{eq:simplex-coefficient-inequality}, we obtain
$v_a(u^*)\ge\frac{\kappa_\delta-\zeta}{C}$ for every $a$.
The initial choices give $\kappa_\delta-\zeta\ge1-\eta$.
Fix $z_a=\frac{u_a^*+\frac{\zeta}{n_0}}{x_a}$ before sampling the construction.
Equation~\eqref{eq:arbitrary-coefficient-bound}, together with $v_a(u^*)\ge\frac{\kappa_\delta-\zeta}{C}\ge\frac{1-\eta}{C}$ for every $a$, gives~\eqref{eq:rounding-marginals}.
For the fixed vector $z_a=\frac{u_a^*+\frac{\zeta}{n_0}}{x_a}$, all constants and selection functions are chosen before the Poisson and uniform variables in \textbf{Definition~\ref{def:local-random-matching}} are drawn. The vector $z$ may depend on the entire original input, but is kept unchanged in every experiment and in every resampling comparison.
This completes the proof.
\end{proof}

\subsection{Changing the random variables of one member}

We use the construction in the proof of \textbf{Proposition~\ref{prop:rounding}} to control how a matching changes when the random variables attached to one member are replaced. Throughout \textbf{Lemma~\ref{lem:rounding-resampling}}, the input $(\mathcal A,x)$, the chosen coefficient vector $z$, all constants, and all selection functions remain fixed.

\begin{lemma}\label{lem:rounding-resampling}
Use the construction of \textbf{Proposition~\ref{prop:rounding}} for one fixed index $d$.
Replace the entire group $(K_c,U_c)$ belonging to a single member $c$ by an independent copy, leaving every other group unchanged, and let $I_d'$ be the resulting matching. There are positive integers $L_0,h$, depending only on $r,\eta$, such that $|I_d\mathbin{\triangle}I_d'|\le L_0$, where $I_d\mathbin{\triangle}I_d'=(I_d\setminus I_d')\cup(I_d'\setminus I_d)$ is the symmetric difference.
For every $b\in I_d\mathbin{\triangle}I_d'$, there is a simple path from $b$ to $c$ of length at most $h$ in $U(\mathcal A)$ such that $K_a>0$ for every member $a\ne c$ on the path in the original experiment.
If $K_c'$ is the resampled value of $K_c$, then $K_c>0$ or $K_c'>0$, unless $b=c$ and the path has length zero.
\end{lemma}

\begin{proof}
Let $K'$ be the vector obtained by replacing $K_c$ with $K_c'$. Among the Poisson variables, only $K_c$ can change; the auxiliary variables attached to $c$ may change as well.
The sum $T_v=\sum_{a\ni v}K_a$ can therefore change only for $v\in c$.
For $a\ne c$, the weight $\frac{K_a}{M}$ of $a$ is unchanged whenever $a$ is retained in both samples.
The retention status of $a\ne c$ can change only if $a\cap c\ne\varnothing$, by~\eqref{eq:retained-family}.
At every original vertex there are at most $M$ retained members in each experiment, since $\sum_{a\in\mathcal A_K:a\ni v}K_a\le M$ and every retained member has $K_a\ge1$, with the corresponding inequalities for $K'$.
Consequently the members whose retention status or attached variables change lie in a set of size at most $1+2rM$, consisting of $c$ and the members of $\mathcal A_K\cup\mathcal A_{K'}$ that meet $c$.

Let $G=U(\mathcal A)[\mathcal A_K\cup\mathcal A_{K'}\cup\{c\}]$ be the subgraph of the line graph induced by the members retained in at least one experiment, together with $c$.
The graph $G-c$ has maximum degree at most $2rM$: each member contains at most $r$ vertices, and each vertex belongs to at most $M$ retained members in each experiment. Adding $c$ increases each other degree by at most one, and $c$ has at most $2rM$ neighbors. Hence $\Delta(G)\le2rM+1\le2rM+r$.
By part~\textup{(iv)} of \textbf{Definition~\ref{def:local-random-matching}}, each member whose selection indicator changes lies within $G$-distance $R_*$ of a member whose retention status or attached variables change.
Indeed, if the retention status and variables of every member within $G$-distance $R_*$ of $b$ agree before and after resampling, then $B_K(b,R_*)=B_{K'}(b,R_*)$ and the attached variables agree on the common ball. Part~\textup{(iv)} of \textbf{Definition~\ref{def:local-random-matching}} therefore forces the selection indicators of $b$ to agree.
Writing $B_0=2rM+r$, we may take $L_0=(1+2rM)\sum_{\ell=0}^{R_*}B_0^\ell$: the expression bounds the total number of vertices in the radius-$R_*$ balls in $G$ centered at the at most $1+2rM$ members whose retention status or attached variables can change. The resulting $L_0$ depends only on $r,M,R_*$.

If a member $b\in I_d\mathbin{\triangle}I_d'$ equals $c$, use a path of length zero.
Otherwise a path of length at most $R_*$ in $G$ joins $b$ to a member whose retention status or attached variables change; one additional step joins the endpoint to $c$ if necessary.
For every member $a\ne c$ of $G$, we have $K_a=K_a'>0$.
If changing $K_c$ causes the retention status of another member to change, then $K_c>0$ or $K_c'>0$.
If only auxiliary variables at $c$ cause a change, then $c$ must have been retained in both experiments, so $K_c=K_c'>0$.
Delete repeated vertices from the path to obtain a simple path of length at most $R_*+1$.
Thus $h=R_*+1$ satisfies \textbf{Lemma~\ref{lem:rounding-resampling}}.
\end{proof}

\section{Partitioning each edge by codegrees}\label{s:partition}

Throughout \textbf{Sections~\ref{s:partition}--\ref{s:auxiliary}}, let $k\ge r\ge2$ be fixed integers, let $\rho,\xi>0$, and let $H$ be a $k$-uniform, $r$-simple multihypergraph whose maximum degree $\Delta=\Delta(H)$ is sufficiently large in terms of $k,r,\rho,\xi$. We prove \textbf{Lemma~\ref{lem:auxiliary}}. Replacing $\xi$ by $\min\{\xi,1\}$, if necessary, we assume that $0<\xi\le1$. Put $C=r-1+\frac{1}{r}$. In \textbf{Section~\ref{s:partition}}, we choose constants and partition the edges of $H$; no random variables are sampled until \textbf{Section~\ref{s:joint}}.

For $0<\theta<\frac{1}{10(k+1)}$, define
\begin{equation}\label{eq:theta-constants}
 A_\theta=\frac{C(1+\theta)}{(1-k\theta)(1-\theta)},\qquad
 B_\theta=(1-k\theta)(1-\theta)-\theta.
\end{equation}

By~\eqref{eq:theta-constants}, $B_\theta>0$ for $0<\theta<\frac{1}{10(k+1)}$. Since $A_\theta\to C$ and $B_\theta\to1$ as $\theta\to0$, we fix $\theta$ sufficiently small that
\begin{equation}\label{eq:theta-choice}
 \frac{A_\theta+2\theta}{B_\theta}+\theta<C+\frac{\rho}{2}.
\end{equation}

Fix the integers $P,m,M,R_*$ supplied by \textbf{Proposition~\ref{prop:rounding}} for $r$ and $\eta=\theta^2$. The integers $P,m,M,R_*$ are independent of the codegree thresholds $\delta_i$ in~\eqref{eq:threshold-choice}. \textbf{Lemma~\ref{lem:weighted-variance}} will give a constant $K>0$ depending only on $k,r,\theta,P,m,M,R_*$. The proof of \textbf{Lemma~\ref{lem:weighted-variance}} does not use the values of $\delta_i$ or $\nu_i$.

Put $m_0=1+\binom{k}{2}$ and $\delta_1=1$. We choose $\delta_2,\ldots,\delta_{m_0+1}$ and $\nu_1,\ldots,\nu_{m_0}$ recursively. At step $i\in[m_0]$, the number $\delta_i>0$ has already been chosen, and we proceed as follows.
\begin{enumerate}[label=(\arabic*)]
\item Let $\Gamma_i$ be the simple graph on $V(H)$ in which $uv$ is an edge if and only if $u\ne v$ and $d_H(u,v)\ge\delta_i\Delta$. For every $v\in V(H)$,
\begin{equation}\label{eq:codegree-neighborhood}
 \delta_i\Delta\,d_{\Gamma_i}(v)
 \le\sum_{u\ne v}d_H(u,v)
 =(k-1)d_H(v)
 \le(k-1)\Delta.
\end{equation}

By~\eqref{eq:codegree-neighborhood}, $\Delta(\Gamma_i)\le L_i$, where $L_i=\max\{1,\lceil\frac{k-1}{\delta_i}\rceil\}$. Set $T_i=\sum_{j=0}^{k-1}L_i^j$ and $F_i=2^{T_i}$. If $v\in U\subseteq V(H)$, $|U|\le k$, and $\Gamma_i[U]$ is connected, then $U$ lies in the ball of radius $k-1$ about $v$ in $\Gamma_i$. The radius-$(k-1)$ ball about $v$ has at most $T_i$ vertices. Consequently, there are at most $F_i$ sets $U$ containing $v$ for which $|U|\le k$ and $\Gamma_i[U]$ is connected.

\item Choose $\nu_i>0$ with $\nu_i\le\min\{\frac{1}{4},\frac{\theta}{m_0F_i}\}$. Define
\begin{equation}\label{eq:level-probabilities}
 \begin{aligned}
 q_i&=(1-k\theta)\left(\frac{\theta\nu_i}{1+\theta}\right)^k,\\
 m_i^*&=(1-k\theta)\min\left\{\theta,\frac{1-\theta}{C}\right\}
                 \frac{\nu_i}{1+\theta}.
 \end{aligned}
\end{equation}

Both numbers in~\eqref{eq:level-probabilities} are positive.
\item Choose $0<\delta_{i+1}<\frac{\delta_i}{2}$ sufficiently small that
\begin{equation}\label{eq:threshold-choice}
 \frac{K\delta_{i+1}}{q_i^2}
 \le\min\left\{\theta^2m_i^*,\frac{\theta^3q_i}{k}\right\},
 \qquad
 \frac{\delta_{i+1}}{B_\theta q_i}\le\frac{\xi}{4}.
\end{equation}

The quantities $K,q_i,m_i^*,B_\theta,\theta,\xi$ have already been fixed at step $i$, so a positive choice of $\delta_{i+1}$ exists.
\end{enumerate}

The bounds $L_i,T_i,F_i$ use only $k$ and $\delta_i$, even though $\Gamma_i$ is defined using $H$. Thus all the chosen constants depend only on $k,r,\rho,\xi$, and none depend on $\Delta$, $|V(H)|$, or $|E(H)|$.

For each edge $e\in E(H)$, let $i(e)$ be the least $i\in[m_0]$ such that no pair of distinct vertices $u,v\in e$ satisfies $\delta_{i+1}\le \frac{d_H(u,v)}{\Delta}<\delta_i$. The index $i(e)$ exists because the $m_0$ intervals $[\delta_{i+1},\delta_i)$ are pairwise disjoint and $e$ contains only $\binom{k}{2}<m_0$ pairs of vertices. We call $i(e)$ the \emph{level} of $e$.

If $i=i(e)$, partition $e$ into the vertex sets of the connected components of $\Gamma_i[e]$. In the construction, a \emph{part} is a pair $a=(i,U)$ for which $i=i(e)$ and $U$ is the vertex set of a connected component of $\Gamma_i[e]$ for some $e\in E(H)$. We call $i$ the \emph{level} of $a$ and write $U(a)=U$ and $|a|=|U|$. Let $B(e)$ be the set of parts obtained from $e$. The pair $(i,U)$ is counted as one part even when $(i,U)\in B(e)\cap B(f)$ for distinct edges $e,f$; pairs with different levels are different parts. Let $\mathcal A$ be the family of sets $U(a)$, labeled by the parts $a=(i,U)$. For $a\in\mathcal A$, define
\begin{equation}\label{eq:part-weights}
 \begin{aligned}
 d_a&=|\{e\in E(H):a\in B(e)\}|,\qquad w_a=\frac{d_a}{\Delta},\\
 v_a&=\frac{w_a+\nu_i}{1+\theta}\quad\text{if $a$ has level $i$}.
 \end{aligned}
\end{equation}

Parts $a,b\in\mathcal A$ intersect if $U(a)\cap U(b)\ne\varnothing$. For each $e$, the sets $U(a)$, $a\in B(e)$, are nonempty, pairwise disjoint, and have union $e$; in particular $1\le|B(e)|\le k$.

We record four properties of the parts $a\in\mathcal A$ and the weights in~\eqref{eq:part-weights}.
\begin{enumerate}[label=(\roman*)]
\item\label{part:weights} For every $u\in V(H)$, $\sum_{a:u\in U(a)}w_a=\frac{d_H(u)}{\Delta}\le1$. Indeed, every edge containing $u$ contributes to exactly one part containing $u$.
\item\label{part:large} If $|a|>r$, then $d_a=1$. Otherwise two distinct edges would both contain $U(a)$ and hence would intersect in more than $r$ vertices, contrary to $r$-simplicity.
\item\label{part:pair} If distinct parts $a,b$ belong to $B(e)$ and $i=i(e)$, then fewer than $\delta_{i+1}\Delta$ edges $f$ satisfy $a,b\in B(f)$. To prove item~\ref{part:pair}, choose $u\in U(a)$ and $v\in U(b)$. The vertices $u,v$ are in different components of $\Gamma_i[e]$, so $d_H(u,v)<\delta_i\Delta$. The definition of $i(e)$ excludes $\frac{d_H(u,v)}{\Delta}\in[\delta_{i+1},\delta_i)$, and therefore $d_H(u,v)<\delta_{i+1}\Delta$. Every edge $f$ with $a,b\in B(f)$ contains $u$ and $v$.
\item\label{part:positive} For every $u\in V(H)$, $\sum_{a:u\in U(a)}v_a\le1$, and every part $a$ of level $i$ satisfies $v_a\ge\frac{\nu_i}{1+\theta}$. To bound $\sum_{a:u\in U(a)}v_a$, note that $\Gamma_i[U(a)]$ is connected and $|U(a)|\le k$ for each part $a$ of level $i$. Thus at most $F_i$ parts of level $i$ contain $u$. By item~\ref{part:weights} and the choice of $\nu_i$,
\begin{equation}\label{eq:positive-weight-sum}
 \sum_{a:u\in U(a)}v_a
 \le\frac{1+\sum_{i=1}^{m_0}F_i\nu_i}{1+\theta}
 \le\frac{1+\theta}{1+\theta}=1.
\end{equation}

The lower bound $v_a\ge\frac{\nu_i}{1+\theta}$ follows from~\eqref{eq:part-weights} and $w_a\ge0$.
\end{enumerate}

\section{Selecting disjoint parts and bounding weighted sums}\label{s:joint}

Apply \textbf{Proposition~\ref{prop:rounding}} to the indexed family $\{a\in\mathcal A:|a|\le r\}$ with weights $x_a=v_a$ and error $\eta=\theta^2$. Equation~\eqref{eq:positive-weight-sum} verifies the fractional matching hypothesis of \textbf{Proposition~\ref{prop:rounding}}. For $d\in[P]$, let $I_d$ be the resulting random matching, and write $p_{a,d}=\PP(a\in I_d)$ when $|a|\le r$. Choose the coefficients supplied by \textbf{Proposition~\ref{prop:rounding}} for the fixed family $\{a\in\mathcal A:|a|\le r\}$ with weights $v_a$, and use the same coefficients in every experiment.

For each part $a$ of size at most $r$, define $\mathcal D_a=\{d\in[P]:p_{a,d}\ge\frac{(1-\theta)v_a}{C}\}$. For $|a|>r$, define $\mathcal D_a=[P]$. The marginal bounds in \textbf{Proposition~\ref{prop:rounding}} give, for $|a|\le r$,
\begin{equation}\label{eq:good-index-count}
 0\le1-\frac{Cp_{a,d}}{v_a},\qquad
 \frac{1}{P}\sum_{d=1}^P\left(1-\frac{Cp_{a,d}}{v_a}\right)\le\theta^2.
\end{equation}

For $d\notin\mathcal D_a$, we have $1-\frac{Cp_{a,d}}{v_a}>\theta$. By~\eqref{eq:good-index-count}, all terms $1-\frac{Cp_{a,d}}{v_a}$ are nonnegative and their average over $d\in[P]$ is at most $\theta^2$, so $|[P]\setminus\mathcal D_a|\le\theta P$. Put $\mathcal F_d=\{f\in E(H):d\in\mathcal D_a\text{ for every }a\in B(f)\}$. Since $|B(f)|\le k$, at most $k\theta P$ indices are excluded by one of the parts of a specified edge $f$. Thus each $f$ belongs to $\mathcal F_d$ for at least $(1-k\theta)P$ values of $d$.

Fix an index $d$. Independently of the variables producing $I_d$, and independently over all $a\in\mathcal A$, choose Bernoulli variables $T_a$ with $\PP(T_a=1)=\theta v_a$. Define a new random family $S_d$ by the following two rules.
\begin{enumerate}[label=(\roman*)]
\item\label{rule:isolated-proposals} Include a part $a$ with $T_a=1$ if $T_b=0$ for every other part $b$ intersecting $a$.
\item\label{rule:retained-matching} Include a part $a\in I_d$ if $T_b=0$ for every part $b$ intersecting $a$, including $b=a$.
\end{enumerate}

The parts included by rule~\ref{rule:isolated-proposals} are pairwise disjoint. Rule~\ref{rule:retained-matching} retains a subfamily of the matching $I_d$ and excludes every member of $I_d$ intersecting a part included by rule~\ref{rule:isolated-proposals}. Consequently $S_d$ is a matching on $\mathcal A$.

For each part $a$ and edge $f$, put $m_{a,d}=\PP(a\in S_d)$ and $q_{f,d}=\PP(B(f)\subseteq S_d)$. The probabilities $m_{a,d},q_{f,d}$ are computed with the family $\mathcal A$, the weights $v_a$, the coefficients from \textbf{Proposition~\ref{prop:rounding}}, and the index $d$ fixed. In particular, $p_{a,d}$, $m_{a,d}$, $q_{f,d}$, and $\mathcal F_d$ do not depend on any sampled outcome and remain unchanged when one group of random variables is replaced by an independent copy.

The inequalities in item~\ref{part:positive} imply $\sum_{b:U(b)\cap U(a)\ne\varnothing}v_b\le |a|\le k$. The union bound, and the independence of $I_d$ from the variables $T_b$, therefore give
\begin{equation}\label{eq:part-marginal-lower}
 \begin{aligned}
 m_{a,d}&\ge\frac{(1-k\theta)(1-\theta)v_a}{C}
       &&\text{if $|a|\le r$ and $d\in\mathcal D_a$},\\
 m_{a,d}&\ge(1-k\theta)\theta v_a
       &&\text{if $|a|>r$}.
 \end{aligned}
\end{equation}

To obtain the bound for $|a|>r$ in~\eqref{eq:part-marginal-lower}, require $T_a=1$ and $T_b=0$ for every $b\ne a$ intersecting $a$, and apply rule~\ref{rule:isolated-proposals}.

If $f$ has level $i$, require $T_a=1$ for every $a\in B(f)$ and $T_b=0$ for every $b\notin B(f)$ with $U(b)\cap f\ne\varnothing$. The parts in $B(f)$ are pairwise disjoint. The variables $T_a$, $a\in B(f)$, are independent of the variables $T_b$, $b\notin B(f)$, and $\sum_{b:U(b)\cap f\ne\varnothing}v_b\le k$ by item~\ref{part:positive}. Hence
\begin{equation}\label{eq:joint-probability-lower}
 q_{f,d}\ge(1-k\theta)\prod_{a\in B(f)}\theta v_a
 \ge(1-k\theta)\left(\frac{\theta\nu_i}{1+\theta}\right)^k=q_i>0.
\end{equation}

The product bound in~\eqref{eq:joint-probability-lower} uses $|B(f)|\le k$ and $0<\frac{\theta\nu_i}{1+\theta}<1$; the equality to $q_i$ follows from~\eqref{eq:level-probabilities}.

If $a$ has level $i$ and $d\in\mathcal D_a$, equations~\eqref{eq:part-weights} and \eqref{eq:part-marginal-lower} give $m_{a,d}\ge m_i^*$. If $|a|\le r$, equations~\eqref{eq:part-weights} and \eqref{eq:part-marginal-lower} also give $\frac{d_a}{m_{a,d}}\le A_\theta \Delta$. For $|a|>r$, item~\ref{part:large} gives $d_a=1$, whereas \eqref{eq:part-marginal-lower} gives $m_{a,d}\ge m_i^*$. Combining the cases $|a|\le r$ and $|a|>r$, we obtain, for $\Delta\ge\max_{i\in[m_0]}\frac{1}{A_\theta m_i^*}$,
\begin{equation}\label{eq:degree-over-marginal}
 m_{a,d}\ge m_i^*,\qquad \frac{d_a}{m_{a,d}}\le A_\theta \Delta.
\end{equation}

For each part $a$, define
\begin{equation}\label{eq:weighted-sum}
 W_{a,d}=\sum_{\substack{f\in\mathcal F_d\\a\in B(f)}}
          \frac{\mathbf1_{\{B(f)\setminus\{a\}\subseteq S_d\}}}{q_{f,d}}.
\end{equation}

If $a$ has level $i$, then every edge in the sum~\eqref{eq:weighted-sum} has level $i$ and $q_{f,d}\ge q_i$ by~\eqref{eq:joint-probability-lower}. Thus $0\le W_{a,d}\le \frac{d_a}{q_i}\le \frac{\Delta}{q_i}$, and hence $\EE W_{a,d}^2<\infty$.
If $m_{a,d}>0$, then $\PP(B(f)\setminus\{a\}\subseteq S_d\mid a\in S_d)=\frac{q_{f,d}}{m_{a,d}}$ for every summand of~\eqref{eq:weighted-sum}, since $a\in B(f)$. For each summand of~\eqref{eq:weighted-sum}, divide the conditional probability $\frac{q_{f,d}}{m_{a,d}}$ by $q_{f,d}$. Summing over $f\in\mathcal F_d$ with $a\in B(f)$ gives
\begin{equation}\label{eq:weighted-conditional-mean}
 \EE(W_{a,d}\mid a\in S_d)
 =\frac{|\{f\in\mathcal F_d:a\in B(f)\}|}{m_{a,d}}.
\end{equation}

\begin{lemma}\label{lem:weighted-variance}
There exists $K>0$, depending only on $k,r,\theta$ and the constants in \textbf{Proposition~\ref{prop:rounding}}, such that for every $d\in[P]$ and every part $a$ of level $i$,
\begin{equation}\label{eq:weighted-variance}
 \Var(W_{a,d})\le \frac{K\delta_{i+1}\Delta^2}{q_i^2}.
\end{equation}

In particular, $K$ is independent of $\delta_1,\ldots,\delta_{m_0+1}$ and $\nu_1,\ldots,\nu_{m_0}$.
\end{lemma}

\begin{proof}
Fix $d$. For each part $b$, let $Y_b$ be the random vector consisting of all variables assigned to $b$. If $|b|\le r$, then $Y_b$ consists of the Poisson variable $K_b$ and the uniform variables in \textbf{Definition~\ref{def:local-random-matching}}, together with $T_b$; if $|b|>r$, then $Y_b$ consists only of $T_b$. The vectors $Y_b$, $b\in\mathcal A$, are mutually independent. In this proof, call a part $b$ \emph{active} if either $|b|\le r$ and $K_b>0$, or $T_b=1$. The events that distinct parts are active are independent. With $\kappa=\max\{1,m+\theta\}$, each part $b$ is active with probability at most $\kappa v_b$.

Let $U$ be the line graph on $\mathcal A$, and let $N_U[b]$ denote the closed neighborhood of $b$. By item~\ref{part:positive},
\begin{equation}\label{eq:conflict-weight-sum}
 \sum_{z\in N_U[b]}v_z
 \le\sum_{u\in U(b)}\sum_{z:u\in U(z)}v_z\le k.
\end{equation}

For a part $z$, replace $Y_z$ by an independent copy and keep $Y_b$ unchanged for every $b\ne z$. Let $D_z$ be the set of parts whose membership in $S_d$ changes. We first establish two properties with constants $L,h$ depending only on $k,r,\theta$:
\begin{enumerate}[label=(\roman*)]
\item\label{influence:number} $|D_z|\le L$ for every outcome and every $z$;
\item\label{influence:path} if $b\in D_z$ and $b\ne z$, there is a simple path of length at most $h$ from $b$ to $z$ in $U$. Every internal vertex on the path from $b$ to $z$ is active in the original experiment, and $z$ is active in at least one of the original and resampled experiments.
\end{enumerate}

If $|z|>r$, then $Y_z$ contains only $T_z$, so replacing $Y_z$ cannot change $I_d$. If $|z|\le r$, apply \textbf{Lemma~\ref{lem:rounding-resampling}} with $x_b=v_b$ on $\{b\in\mathcal A:|b|\le r\}$, using the coefficients fixed when we applied \textbf{Proposition~\ref{prop:rounding}}. The conclusion of \textbf{Lemma~\ref{lem:rounding-resampling}} bounds the number of changes in membership in $I_d$ and joins each changed part to $z$ by a path whose internal parts have positive Poisson variables. A part with a positive Poisson variable is active. \textbf{Lemma~\ref{lem:rounding-resampling}} also accounts for changes to the probabilities used in the sequential inclusion decisions, because the inclusion probabilities are functions of the local variables specified in the proof of \textbf{Proposition~\ref{prop:rounding}}. With all $T_b$ fixed, rule~\ref{rule:retained-matching} retains exactly the members of $I_d$ meeting no part $b$ with $T_b=1$, and rule~\ref{rule:isolated-proposals} gives the same selected parts in both experiments. Consequently, changing $I_d$ with all $T_b$ fixed can change membership in $S_d$ only at a part whose membership in $I_d$ changes.

Now change $T_z$ with $I_d$ fixed. Every changed part intersects $z$. The parts included by rule~\ref{rule:isolated-proposals} before resampling form a matching, and the parts included by rule~\ref{rule:isolated-proposals} after resampling form a matching. In each matching, at most $|z|\le k$ members intersect $z$. At most $k$ members of $I_d$ intersect $z$, so at most $3k$ membership decisions change. A changed part $b\ne z$ is active, and the one-edge path $bz$ satisfies item~\ref{influence:path}. Replacing $Y_z$ can be performed by first changing the variables producing $I_d$ with all $T_b$ fixed, and then changing $T_z$ with $I_d$ fixed. Combining the bounds for changes to $I_d$ and $T_z$ proves items~\ref{influence:number} and~\ref{influence:path}, after increasing $L$ and $h$. The constants $L,h$ depend only on $k,r,\theta$ and not on the thresholds $\delta_i,\nu_i$ or on a positive lower bound for $v_b$.

For each fixed $b$, item~\ref{influence:path} and independence of the vectors $Y_z$ give
\begin{equation}\label{eq:summed-influence}
 \sum_z\PP(b\in D_z)
 \le K_2:=1+2\sum_{\ell=1}^h(\kappa k)^\ell.
\end{equation}

To prove~\eqref{eq:summed-influence}, bound the term $z=b$ by one. For $z\ne b$, enumerate simple paths $b=b_0,b_1,\ldots,b_\ell=z$ of length $\ell\le h$ as in item~\ref{influence:path}. The probability that $b_1,\ldots,b_{\ell-1}$ are active in the original experiment and $z$ is active before or after resampling is at most $2\kappa^\ell\prod_{j=1}^\ell v_{b_j}$. The factor two accounts for activity of $z$ before or after resampling. Since the vertices on the path are distinct, the vectors $Y_{b_1},\ldots,Y_{b_\ell}$ are independent, and the independent copy replacing $Y_z$ is independent of the internal vectors $Y_{b_1},\ldots,Y_{b_{\ell-1}}$. The sum of $\kappa^\ell\prod_{j=1}^\ell v_{b_j}$ over simple paths is at most the corresponding sum over all walks of length $\ell$. Repeated use of~\eqref{eq:conflict-weight-sum}, starting at the final step of a walk, bounds the sum over walks by $(\kappa k)^\ell$. The union bound over $\ell\in[h]$ and the term $z=b$ prove~\eqref{eq:summed-influence}.

Set $A_a=0$, and for $b\ne a$ set
\begin{equation}\label{eq:influence-coefficients}
 A_b=\sum_{\substack{f\in\mathcal F_d\\a,b\in B(f)}}\frac{1}{q_{f,d}}.
\end{equation}

Every $f$ appearing in~\eqref{eq:influence-coefficients} has the level $i$ of $a$. Item~\ref{part:pair}, equation~\eqref{eq:joint-probability-lower}, and $d_a\le \Delta$ imply
\begin{equation}\label{eq:influence-coefficient-bounds}
 \begin{aligned}
 \max_b A_b&\le\frac{\delta_{i+1}\Delta}{q_i},\\
 \sum_b A_b
 &=\sum_{\substack{f\in\mathcal F_d\\a\in B(f)}}\frac{|B(f)|-1}{q_{f,d}}
 \le\frac{k\Delta}{q_i}.
 \end{aligned}
\end{equation}

Write $W_{a,d}^{(z)}$ for the value of~\eqref{eq:weighted-sum} after replacing $Y_z$ by an independent copy. A summand of \eqref{eq:weighted-sum} changes only if some $b\in B(f)\setminus\{a\}$ belongs to $D_z$. Therefore
\begin{equation}\label{eq:weighted-change}
 \begin{aligned}
 |W_{a,d}-W_{a,d}^{(z)}|&\le\sum_{b\in D_z}A_b,\\
 (W_{a,d}-W_{a,d}^{(z)})^2&\le L\sum_b A_b^2\mathbf1_{\{b\in D_z\}}.
 \end{aligned}
\end{equation}

Applying \textbf{Lemma~\ref{lem:efron-stein}} to the independent vectors $Y_z$, and then using~\eqref{eq:weighted-change}, \eqref{eq:summed-influence}, and \eqref{eq:influence-coefficient-bounds}, we obtain
\begin{equation}\label{eq:variance-derivation}
 \begin{aligned}
 \Var(W_{a,d})
 &\le\frac{L}{2}\sum_b A_b^2\sum_z\PP(b\in D_z)
 \le\frac{LK_2}{2}\sum_b A_b^2\\
 &\le\frac{LK_2}{2}(\max_b A_b)\sum_b A_b
 \le\frac{LK_2k}{2}\,\frac{\delta_{i+1}\Delta^2}{q_i^2}.
 \end{aligned}
\end{equation}

Taking $K=\frac{LK_2k}{2}$ in~\eqref{eq:variance-derivation} proves~\eqref{eq:weighted-variance}. Since $L,h,\kappa,k$ are independent of $\delta_i,\nu_i$, the formula for $K_2$ in~\eqref{eq:summed-influence} also proves the asserted independence of $K$ from $\delta_i,\nu_i$.
\end{proof}

We now use the choices in \eqref{eq:threshold-choice}. Suppose that $a$ has level $i$ and $d\in\mathcal D_a$, and put $E_a=\{a\in S_d\}$. For any square-integrable random variable $X$ and event $E$ of positive probability, the Cauchy--Schwarz inequality gives $|\EE(X\mid E)-\EE X|\le\sqrt{\frac{\Var(X)}{\PP(E)}}$, by writing the difference as $\frac{\EE((X-\EE X)\mathbf1_E)}{\PP(E)}$. Consequently, equations~\eqref{eq:degree-over-marginal}, \eqref{eq:weighted-conditional-mean}, \eqref{eq:weighted-variance}, and \eqref{eq:threshold-choice} imply
\begin{equation}\label{eq:weighted-mean-bound}
 |\EE W_{a,d}-\EE(W_{a,d}\mid E_a)|
 \le\Delta\sqrt{\frac{K\delta_{i+1}}{q_i^2m_i^*}}\le\theta\Delta,
 \qquad \EE W_{a,d}\le(A_\theta+\theta)\Delta.
\end{equation}

For $f\in\mathcal F_d$ with $a\in B(f)$, the mean bound~\eqref{eq:weighted-mean-bound}, Chebyshev's inequality, and~\eqref{eq:joint-probability-lower} give
\begin{equation}\label{eq:conditional-overflow}
 \begin{aligned}
 &\PP\bigl(W_{a,d}>(A_\theta+2\theta)\Delta\mid B(f)\subseteq S_d\bigr)\\
 &\quad\le\frac{\PP(W_{a,d}-\EE W_{a,d}>\theta \Delta)}{q_{f,d}}
 \le\frac{K\delta_{i+1}}{q_i^3\theta^2}
 \le\frac{\theta}{k}.
 \end{aligned}
\end{equation}

\begin{definition}\label{def:edge-retention}
In an experiment with index $d$, an edge $f$ is \emph{retained} if $f\in\mathcal F_d$, $B(f)\subseteq S_d$, and $W_{a,d}\le(A_\theta+2\theta)\Delta$ for every $a\in B(f)$.
\end{definition}

By \textbf{Definition~\ref{def:edge-retention}}, the union bound over the $|B(f)|\le k$ parts in~\eqref{eq:conditional-overflow} gives
\begin{equation}\label{eq:edge-retention}
 \PP(f\text{ is retained in an experiment with index }d)
 \ge(1-\theta)q_{f,d}\quad(f\in\mathcal F_d).
\end{equation}

\section{Partition into hypergraphs of small codegree}\label{s:auxiliary}

\begin{proof}[Proof of \textbf{Lemma~\ref{lem:auxiliary}}]
Use the constants, parts, and experiments from \textbf{Sections~\ref{s:partition} and~\ref{s:joint}}. Let $s=P\lfloor\frac{\sqrt \Delta}{P}\rfloor$. Since $s\ge\sqrt\Delta-P$, the assumption $\Delta\ge4P^2$ gives $\frac{\sqrt\Delta}{2}\le s\le\sqrt\Delta$ and $s>0$. List each index $d\in[P]$ exactly $\frac{s}{P}$ times to obtain a fixed list $d_1,\ldots,d_s$. In experiment $j$, apply the construction of $S_d$ with $d=d_j$, using an independent copy of the full family of Poisson, uniform, and Bernoulli variables specified in \textbf{Section~\ref{s:joint}}. The collections of variables used by different experiments are mutually independent, even when the indices $d_j$ coincide. Denote the selected matching in experiment $j$ by $S^{(j)}$, and the value of~\eqref{eq:weighted-sum} in experiment $j$ by $W_a^{(j)}$. The fixed probabilities $q_{f,d_j}$ are the same for experiments with the same index. Let $R_{f,j}$ be the indicator that $f$ is retained in experiment $j$, and put
\begin{equation}\label{eq:assignment-normalizer}
 Z_f=\sum_{j=1}^s\frac{R_{f,j}}{q_{f,d_j}}.
\end{equation}

Each $f$ belongs to $\mathcal F_d$ for at least $(1-k\theta)P$ indices. By~\eqref{eq:edge-retention}, $\EE Z_f\ge(1-k\theta)(1-\theta)s=(B_\theta+\theta)s$. If $f$ has level $i$, the summands in \eqref{eq:assignment-normalizer} are independent and lie in $[0,\frac{1}{q_i}]$, by \eqref{eq:joint-probability-lower}. \textbf{Lemma~\ref{lem:hoeffding}} therefore gives
\begin{equation}\label{eq:normalizer-tail}
 \PP(Z_f<B_\theta s)\le\exp(-2\theta^2q_i^2s).
\end{equation}

We apply \textbf{Lemma~\ref{lem:lll}} to the events $\{Z_f<B_\theta s\}$ in order to obtain $Z_f\ge B_\theta s$ for every $f\in E(H)$. To bound dependencies among the events $\{Z_f<B_\theta s\}$, define a graph $G$ with vertex set $\mathcal A$: distinct $a,b\in\mathcal A$ are adjacent if $U(a)\cap U(b)\ne\varnothing$ or $a,b\in B(g)$ for some $g\in E(H)$. For every $u\in V(H)$, the number of parts containing $u$ is at most $\Delta$, by item~\ref{part:weights} and $d_a\ge1$. Thus a part $a$ has at most $k\Delta$ neighbors $b$ with $U(a)\cap U(b)\ne\varnothing$. By~\eqref{eq:part-weights}, $a$ belongs to $B(g)$ for $d_a\le\Delta$ edges $g$, each with at most $k-1$ other parts. Hence $\Delta(G)\le2k\Delta$.

Membership of a part $b$ in $S_d$ is determined by the random variables assigned to parts at distance at most $R_0:=R_*+2$ from $b$ in $G$. The integer $R_0$ is independent of $\Delta$ and $H$. If $|b|>r$, then $b\notin I_d$. If $|b|\le r$, the condition $b\in\mathcal A_K$ in~\eqref{eq:retained-family} depends only on the Poisson variables of parts intersecting $b$. When $b\notin\mathcal A_K$, we have $b\notin I_d$ by \textbf{Definition~\ref{def:local-random-matching}}(iii). When $b\in\mathcal A_K$, the ball $B_{U(\mathcal A_K)}(b,R_*)$ is contained in $B_G(b,R_*)$. For each part in $B_G(b,R_*)$, membership in $\mathcal A_K$ is determined by the variables assigned to parts in $B_G(b,R_*+1)$. Thus \textbf{Definition~\ref{def:local-random-matching}}(iv) implies that membership of $b$ in $I_d$ is determined by variables assigned to parts in $B_G(b,R_*+1)$. Rules~\ref{rule:isolated-proposals} and~\ref{rule:retained-matching} additionally depend only on $T_c$ for parts $c$ intersecting $b$, so membership in $S_d$ is determined by variables within distance $R_0=R_*+2$ in $G$.

To decide whether $f$ is retained under \textbf{Definition~\ref{def:edge-retention}}, we must determine whether $a\in S_d$ for every $a\in B(f)$ and evaluate $W_{a,d}$ for every $a\in B(f)$. Every part appearing in a summand of $W_{a,d}$ belongs to $B(g)$ for an edge $g$ with $a\in B(g)$, and is consequently equal or adjacent to $a$ in $G$. By~\eqref{eq:weighted-sum} and \eqref{eq:assignment-normalizer}, the event $\{Z_f<B_\theta s\}$ is therefore determined by the variables from all $s$ experiments assigned to parts within distance $R_0+1$ of $B(f)$ in $G$. The coefficients from \textbf{Proposition~\ref{prop:rounding}}, the probabilities $q_{g,d}$, and the families $\mathcal F_d$ are fixed before sampling the experiments. In particular, the event $\{Z_f<B_\theta s\}$ does not depend on any additional random variables used to determine the coefficients, $q_{g,d}$, or $\mathcal F_d$.

Put $R_1=R_0+1$. If $\{Z_f<B_\theta s\}$ and $\{Z_g<B_\theta s\}$ depend on a common random variable, some part in $B(f)$ is at distance at most $2R_1$ from some part in $B(g)$ in $G$. The union of the radius-$2R_1$ balls centered at the parts of $B(f)$ contains at most $k\sum_{\ell=0}^{2R_1}(2k\Delta)^\ell$ parts. Each part belongs to at most $\Delta$ sets $B(g)$. Thus a dependency graph for the events in \eqref{eq:normalizer-tail} has maximum degree at most
\begin{equation}\label{eq:first-dependency-bound}
 b_1=\left\lceil k\Delta\sum_{\ell=0}^{2R_1}(2k\Delta)^\ell\right\rceil.
\end{equation}

Events determined by disjoint subfamilies of the independent random variables are jointly independent, so the graph with degree bound~\eqref{eq:first-dependency-bound} is a dependency graph for the events $\{Z_f<B_\theta s\}$. All thresholds are fixed, and $q_{\min}=\min_{i\in[m_0]}q_i>0$. Since $s\ge\frac{\sqrt \Delta}{2}$, equations~\eqref{eq:normalizer-tail} and \eqref{eq:first-dependency-bound} give $e(b_1+1)\exp(-2\theta^2q_{\min}^2s)\to0$ as $\Delta\to\infty$. By \textbf{Lemma~\ref{lem:lll}}, we can fix experiment outcomes satisfying $Z_f\ge B_\theta s$ for every $f\in E(H)$.

For each edge $f$, independently choose an index $J_f\in[s]$ according to
\begin{equation}\label{eq:assignment-probability}
 \PP(J_f=j)=\frac{R_{f,j}}{q_{f,d_j}Z_f}.
\end{equation}

The probabilities in~\eqref{eq:assignment-probability} are nonnegative and sum to one by~\eqref{eq:assignment-normalizer}; every denominator $q_{f,d_j}Z_f$ is positive by~\eqref{eq:joint-probability-lower} and $Z_f\ge B_\theta s>0$. If $J_f=j$, assign $f$ to $E_j$. For each $j$, form a multihypergraph $H_j$ with vertex set $S^{(j)}$. For every $f\in E_j$, include a separate edge $\varphi_j(f)$ whose vertex set is $B(f)$. In particular, $\varphi_j(f)$ and $\varphi_j(g)$ are distinct edges when $f\ne g$, even if $B(f)=B(g)$. Thus $\varphi_j:E_j\to E(H_j)$ is a bijection, and every edge of $H_j$ is nonempty and has size at most $k$.

We next choose the independent indices $J_f$ so that every $H_j$ satisfies the degree and codegree bounds of \textbf{Lemma~\ref{lem:auxiliary}}. Fix $j$ and $a\in S^{(j)}$. If no edge $f$ retained in experiment $j$ satisfies $a\in B(f)$, then $d_{H_j}(a)=0$. Otherwise, \textbf{Definition~\ref{def:edge-retention}} gives $W_a^{(j)}\le(A_\theta+2\theta)\Delta$. Moreover, $B(f)\subseteq S^{(j)}$ for every edge $f$ retained in experiment $j$. Thus \eqref{eq:assignment-probability} and \eqref{eq:weighted-sum} give
\begin{equation}\label{eq:aux-degree-mean}
 \begin{aligned}
 \EE d_{H_j}(a)
 &=\sum_{\substack{f:a\in B(f)}}\frac{R_{f,j}}{q_{f,d_j}Z_f}
 \\&\le\frac{1}{B_\theta s}\sum_{\substack{f:a\in B(f)}}\frac{R_{f,j}}{q_{f,d_j}}\\
 &\le\frac{W_a^{(j)}}{B_\theta s}
 \le\frac{A_\theta+2\theta}{B_\theta}\,\frac{\Delta}{s}.
 \end{aligned}
\end{equation}

If distinct parts $a,b$ occur together in some $B(f)$, the levels of $a,b$ coincide, say at $i$. Item~\ref{part:pair} bounds $|\{f\in E(H):a,b\in B(f)\}|$ by $\delta_{i+1}\Delta$. Equations~\eqref{eq:joint-probability-lower}, \eqref{eq:assignment-probability}, and \eqref{eq:threshold-choice} consequently give
\begin{equation}\label{eq:aux-codegree-mean}
 \EE d_{H_j}(a,b)
 \le\frac{\delta_{i+1}\Delta}{B_\theta sq_i}
 \le\frac{\xi \Delta}{4s}.
\end{equation}

If $a,b$ never occur together in a set $B(f)$, then $d_{H_j}(a,b)=0$ for every assignment.

Put $a_0=\frac{A_\theta+2\theta}{B_\theta}$ and $d'=\lceil\frac{(a_0+\theta)\Delta}{s}\rceil$. Every degree and codegree in \eqref{eq:aux-degree-mean}--\eqref{eq:aux-codegree-mean} is a sum of independent Bernoulli variables, because the choices $J_f$ are independent over edges. \textbf{Lemma~\ref{lem:bernstein}} implies
\begin{equation}\label{eq:aux-degree-tail}
 \PP(d_{H_j}(a)>d')\le\exp(-\frac{c_1\Delta}{s}),\qquad
 c_1=\frac{\theta^2}{2(a_0+\frac{\theta}{3})}>0.
\end{equation}

To obtain~\eqref{eq:aux-degree-tail}, use~\eqref{eq:aux-degree-mean} to bound $d'-\EE d_{H_j}(a)$ from below by $\frac{\theta\Delta}{s}$. For the codegree $d_{H_j}(a,b)$, equation~\eqref{eq:aux-codegree-mean} gives $\EE d_{H_j}(a,b)\le\frac{\xi\Delta}{4s}$. Applying \textbf{Lemma~\ref{lem:bernstein}} with a gap of $\frac{3\xi\Delta}{4s}$ therefore yields
\begin{equation}\label{eq:aux-codegree-tail}
 \PP(d_{H_j}(a,b)>\frac{\xi \Delta}{s})\le\exp(-\frac{c_2\Delta}{s}),\qquad c_2=\frac{9\xi}{16}>0.
\end{equation}

Each event in~\eqref{eq:aux-degree-tail} or~\eqref{eq:aux-codegree-tail} depends on $J_f$ for at most $\Delta$ edges $f\in E(H)$. A single variable $J_f$ appears, for each $j$, in at most $|B(f)|$ events in~\eqref{eq:aux-degree-tail} and $\binom{|B(f)|}{2}$ events in~\eqref{eq:aux-codegree-tail}. Joining two events whenever both depend on the same $J_f$ therefore gives a dependency graph of maximum degree at most $b_2=\lceil \Delta s(k+\binom{k}{2})\rceil$. With $c=\min\{c_1,c_2\}$, we have $e(b_2+1)\exp(-\frac{c\Delta}{s})\to0$, since $s\le\sqrt \Delta$. A second application of \textbf{Lemma~\ref{lem:lll}} gives choices for which $\Delta(H_j)\le d'$ and $\Delta_2(H_j)\le\frac{\xi \Delta}{s}$ for every $j$.

We now verify the intersection condition of \textbf{Lemma~\ref{lem:auxiliary}} and the bound $\Delta(H_j)\le\frac{(C+\rho)\Delta}{s}$. Suppose that distinct $f,g\in E_j$ meet at an original vertex $u$. There is a unique part $a\in B(f)$ containing $u$ and a unique part $b\in B(g)$ containing $u$. The parts $a,b$ belong to the matching $S^{(j)}$, because $f$ and $g$ are retained in experiment $j$. If $a\ne b$, then $U(a)\cap U(b)=\varnothing$, contradicting $u\in U(a)\cap U(b)$. Hence $a=b$, so $\varphi_j(f)$ and $\varphi_j(g)$ have a common auxiliary vertex. Conversely, if $a\in\varphi_j(f)\cap\varphi_j(g)$, then $\varnothing\ne U(a)\subseteq f\cap g$.

Finally, \eqref{eq:theta-choice} gives $a_0+\theta<C+\frac{\rho}{2}$. Since $\frac{\Delta}{s}\ge\sqrt\Delta$, the assumption $\Delta\ge\frac{4}{\rho^2}$ makes the additional one allowed by the ceiling defining $d'$ at most $\frac{\rho\Delta}{2s}$. Thus $d'\le\frac{(C+\rho)\Delta}{s}$. The families $E_j$ partition $E(H)$ with multiplicity, the maps $\varphi_j$ are bijections, and the bounds $\Delta(H_j)\le\frac{(C+\rho)\Delta}{s}$ and $\Delta_2(H_j)\le\frac{\xi\Delta}{s}$ verify the degree conditions of \textbf{Lemma~\ref{lem:auxiliary}}.
\end{proof}

\section{Proof of \texorpdfstring{\textbf{Theorem~\ref{thm:main}}}{the main theorem}}\label{s:proof-main}

\begin{proof}[Proof of \textbf{Theorem~\ref{thm:main}}]
Let $H$ satisfy the hypotheses of \textbf{Theorem~\ref{thm:main}}, and put $\Delta=\Delta(H)$.
For $k=1$, give the singleton edges at each vertex distinct
colors from $[\Delta]$. Edges at different vertices are disjoint, so the
assigned colors form a proper edge-coloring with $\Delta$ colors. Suppose that $k\ge2$.
If $t=1$, then $\Delta_2(H)\le1$. Let $\xi>0$ and $\Delta_1$ be
supplied by \textbf{Theorem~\ref{thm:small-codegree}} for rank $k$ and
error $\varepsilon$. When
$\Delta\ge\max\{\Delta_1,\frac{1}{\xi}\}$, we have
$\Delta_2(H)\le1\le\xi\Delta$, so
\textbf{Theorem~\ref{thm:small-codegree}} gives~\eqref{eq:main-bound}.

Now let $r=t\ge2$ and $C=r-1+\frac{1}{r}$. Choose $\alpha,\rho>0$
such that $(1+\alpha)(C+\rho)<C+\frac{\varepsilon}{2}$. Let $\xi>0$
and $\Delta_*$ be supplied by \textbf{Theorem~\ref{thm:small-codegree}}
for rank $k$ and error $\alpha$. Apply
\textbf{Lemma~\ref{lem:auxiliary}} with parameters $k,r,\rho,\xi$.
Put $q=\lceil\frac{(C+\rho)\Delta}{s}\rceil$. For each $j\in[s]$,
add to $H_j$ a new vertex $z_j$ and $q$ new edges, each with
underlying set $\{z_j\}$. Call the resulting multihypergraph $J_j$.
By \textbf{Lemma~\ref{lem:auxiliary}}(ii), $\Delta(H_j)\le q$.
The added singleton edges are disjoint from every edge of $H_j$ and contain no pair
of distinct vertices. Since $C>1$, we also have $q\ge\frac{\Delta}{s}$.
By \textbf{Lemma~\ref{lem:auxiliary}}(i),
\begin{equation}\label{eq:apply-kahn}
 \Delta(J_j)=q,\qquad
 \Delta_2(J_j)=\Delta_2(H_j)\le\frac{\xi\Delta}{s}\le\xi q,
 \qquad q\ge(C+\rho)\sqrt{\Delta}.
\end{equation}

The rank of $J_j$ is still at most $k$. For sufficiently large $\Delta$,
the bound $q\ge(C+\rho)\sqrt{\Delta}$ in~\eqref{eq:apply-kahn}
gives $q\ge\Delta_*$.
Thus \textbf{Theorem~\ref{thm:small-codegree}} gives a proper edge-coloring
of $J_j$ with at most $(1+\alpha)q$ colors. Restrict the coloring of
$J_j$ to $H_j$.

Use disjoint color sets for different indices $j$, and give $e\in E_j$
the color of $\varphi_j(e)$. If intersecting original edges belong to the
same $E_j$, their images intersect by
\textbf{Lemma~\ref{lem:auxiliary}}(iii), so the original edges receive different colors.
Edges in different families receive colors from disjoint sets. The resulting
coloring is therefore proper. Using the definition of $q$,
$(1+\alpha)(C+\rho)<C+\frac{\varepsilon}{2}$, and
\textbf{Lemma~\ref{lem:auxiliary}}(i), we obtain
\begin{equation}\label{eq:final-color-bound}
 \begin{aligned}
 \chi'(H)&\le s(1+\alpha)q
 \le(1+\alpha)\bigl((C+\rho)\Delta+s\bigr)\\
 &\le\left(C+\frac{\varepsilon}{2}\right)\Delta
       +(1+\alpha)\sqrt{\Delta}
 \le(C+\varepsilon)\Delta.
 \end{aligned}
\end{equation}

The bound $(1+\alpha)\sqrt{\Delta}\le\frac{\varepsilon\Delta}{2}$ used
in~\eqref{eq:final-color-bound} holds once
$\Delta\ge\left(\frac{2(1+\alpha)}{\varepsilon}\right)^2$.
Thus~\eqref{eq:final-color-bound} proves~\eqref{eq:main-bound}.
\end{proof}

\bigskip
\noindent\textbf{Acknowledgement}

This research was supported in part by the National Natural Science Foundation
of China (Nos. 12471330 and 12571373) and the Youth Student Fundamental
Research Funds of Shandong University.

\par\medskip
\noindent\begin{minipage}{\linewidth}
\setlength{\parindent}{1.5em}
\noindent\textbf{Declaration on the Use of Generative AI}\par\smallskip

We used the gpt-6.1-sol to assist with the exploration of some proof strategies
and calculations and with auxiliary tasks, including checking typographical
and linguistic errors, improving the organization and clarity of the manuscript,
and drawing attention to possible inconsistencies or gaps that were subsequently
examined by the authors. The authors independently verified all statements
and proofs and take full responsibility for the content of this paper.\par
\end{minipage}\par

\begingroup
\small
\setlength{\bibsep}{3pt plus .3pt minus .3pt}
\bibliography{references}
\endgroup
\end{document}